\documentclass[11pt,a4paper]{amsart}
\usepackage{amssymb}
\usepackage{amscd}
\usepackage{amsmath}\usepackage{amsfonts}
\usepackage{pstricks,pst-node}
\usepackage{times,mathrsfs}
\usepackage{cite}

\let\emph=\textbf
\numberwithin{equation}{section}
\newcounter{main}
\theoremstyle{plain}
\newtheorem{Theorem}[main]{Theorem}
\swapnumbers
\newtheorem{theorem}[equation]{Theorem}
\newtheorem*{Corollary}{Corollary}
\newtheorem{corollary}[equation]{Corollary}
\newtheorem{lemma}[equation]{Lemma}

\newtheorem{proposition}[equation]{Proposition}

\newcounter{mycase}
\numberwithin{mycase}{equation}
\def\Case#1{\refstepcounter{mycase}
       \medskip\noindent\textit{Case~\arabic{mycase}}. #1:\space}

\theoremstyle{definition}
\newtheorem{definition}[equation]{Definition}
\newtheorem{example}[equation]{Example}

\theoremstyle{remark}
\newtheorem*{remark}{Remark}

\makeatletter
\def\enumerate{\begingroup\ifnum\@enumdepth>3\@toodeep\else
      \advance\@enumdepth\@ne
      \edef\@enumctr{enum\romannumeral\the\@enumdepth}%
      \topsep\z@\parskip\z@
      \list{\csname label\@enumctr\endcsname}
        {\@nmbrlisttrue\let\@listctr\@enumctr
         \parsep\z@\itemsep\z@\topsep\z@
         \setcounter{\@enumctr}{0}
         \def\makelabel##1{\hss\llap{\rm ##1}}
       }\fi
}
\makeatother

\def\crulefill{\leavevmode\leaders\hrule height 1pt\hfill\kern 0pt}
\long\def\QUERY#1{%
  \leavevmode\newline\noindent$\star\star\star$%
     \thinspace\textsf{Comment/Query}\crulefill\newline%
     \space #1\leavevmode\newline\hbox to 120mm{\crulefill}%
     $\star\star\star$\newline
}

\let\gdom=\triangleright

\let\gedom=\trianglerighteq
\let\bar=\overline

\def\c{\mathbf{c}}

\def\N{\mathbb N}\def\Z{\mathbb Z}

\def\F{\mathbb F}

\newcommand{\Sch}{\ensuremath{\mathscr{S}}}
\newcommand{\h}{\ensuremath{\mathscr{H}}}
\let\H=\h
\newcommand{\Haff}{\ensuremath{\H^{\text{aff}}}}

\newcommand{\Lamb}{\ensuremath{\Lambda^+_{r,n}}}
\def\Sym{\mathfrak S}
\def\Sch{\mathscr S}
\newcommand{\End}{\operatorname{End}}
\newcommand{\Hom}{\operatorname{Hom}}
\newcommand{\Ress}{\operatorname{Res}}
\newcommand{\Res}{\Ress(\Lambda^+_{r,n})}
\newcommand{\la}{\lambda}

\def\nup{\nu_p'}
\newcommand{\bl}{{\boldsymbol \lambda}}
\newcommand{\bmu}{{\boldsymbol \mu}}
\newcommand{\bnu}{{\boldsymbol \nu}}
\newcommand{\brho}{{\boldsymbol \rho}}
\newcommand{\bet}{{\boldsymbol \eta}}
\newcommand{\bxi}{{\boldsymbol \xi}}
\def\O{\mathcal{O}}
\newcommand{\Q}{\mathbf{Q}}
\def\X{\mathbf{X}}

\let\<\langle
\let\>\rangle
\def\({\big(}
\def\){\big)}

\def\Sum{\displaystyle\sum}

\def\map#1#2{\,{:}\,#1\!\longrightarrow\!#2}

\def\abacus#1{%
   \begin{psmatrix}[colsep=2mm,rowsep=1mm,emnode=dot,dotscale=1.3,nodesep=0pt]%
      #1
   \end{psmatrix}
}

{\catcode`\|=\active
  \gdef\set#1{\mathinner{\lbrace\,{\mathcode`\|"8000%
                                   \let|\midvert #1}\,\rbrace}}
}
\def\midvert{\egroup\mid\bgroup}

\DeclareMathOperator{\We}{W_e}
\DeclareMathOperator{\we}{w_e}
\DeclareMathOperator{\Wt}{wt}
\DeclareMathOperator{\ww}{w_e}
\DeclareMathOperator{\winfinity}{w_{\infty}}
\DeclareMathOperator{\leg}{\ell\ell}

\DeclareMathOperator{\rad}{rad}
\DeclareMathOperator{\Ext}{Ext}

\DeclareMathOperator{\res}{res}

\begin{document}
\bibliographystyle{andrew}
\title{Blocks of cyclotomic Hecke algebras}
\author{Sin\'{e}ad Lyle}
\address{School of Mathematics, University of East Anglia, Norwich NR4 7TJ, UK.}
\email{s.lyle@uea.ac.uk}
\author{Andrew Mathas}
\address{School of Mathematics and Statistics F07, University of
Sydney, NSW 2006, Australia.}
\email{a.mathas@maths.usyd.edu.au}


\begin{abstract} This paper classifies the blocks of the cyclotomic Hecke algebras of type $G(r,1,n)$ over an arbitrary field.  Rather than working with the Hecke algebras directly we work instead with the cyclotomic Schur algebras. The advantage of these algebras is that the cyclotomic Jantzen sum formula gives an easy combinatorial characterization of the blocks of the cyclotomic Schur algebras. We obtain an explicit description of the blocks by analyzing the combinatorics of `Jantzen equivalence'.

We remark that a proof of the classification of the blocks of the cyclotomic Hecke algebras was announced in~1999. Unfortunately, Cox has discovered that this previous proof is incomplete.
\end{abstract}
\maketitle


\section{Introduction}
The Ariki--Koike algebras are the cyclotomic
Hecke algebras of type $G(r,1,n)$. These algebras first appeared in
the work of Cherednik~\cite{Cherednik:GT} and they were first
systematically studied by Ariki and Koike~\cite{AK}. Independently,
and at about the same time, Brou\'e and Malle~\cite{BM:cyc}
generalized the definition of Iwahori--Hecke algebras to attach a
Hecke algebra to each complex reflection group. The cyclotomic Hecke
algebras are central to the conjectures of Brou\'e, Malle and
Michel~\cite{Broue:conjectures} which grew out of an attempt to
understand Brou\'e's abelian defect group conjecture for the finite
groups of Lie type.

The Ariki--Koike algebras arise most naturally as `cyclotomic
quotients' of the (extended) affine Hecke algebras of type~$A$.
To make this explicit, let $\F$ be a field and let $\Haff_n$ be the affine
Hecke algebra of type $A_n$. Using the Bernstein presentation,
$\Haff_n$ can be written as a twisted tensor product
$\H_q(\Sym_n)\otimes\F[X_1^\pm,\dots,X_n^\pm]$ of the Iwahori--Hecke
algebra $\H_q(\Sym_n)$ of the symmetric group and the
Laurent polynomial ring $\F[X_1^\pm,\dots,X_n^\pm]$. The Ariki--Koike
algebra is then the quotient algebra
$\H_{r,n}(q,\Q)=\Haff_n/\<(X_1-Q_1)\dots(X_1-Q_r)\>$, where
$\Q=(Q_1,\dots,Q_r)\in(\F^\times)^r$.

As $\H_{r,n}(q,\Q)$ is a quotient of $\Haff_n$, every irreducible
$\H_{r,n}(q,\Q)$--module can be considered as an irreducible
$\Haff_n$--module. Conversely, by quotienting out by the
characteristic polynomial of $X_1$, every irreducible
$\Haff_n$--module is an irreducible module for some Ariki--Koike
algebra. The deep results of Ariki~\cite{Ariki:can} and
Grojnowski~\cite{Groj:control} show that the module categories of the
affine Hecke algebras and the Ariki--Koike algebras are intimately
intertwined. The main result of this paper shows that, combinatorially
at least, the blocks of these algebras are the same.

If $A$ is an algebra then two simple $A$--modules $D$ and
$D'$ belong to same \emph{block} if there exist simple
$A$--modules $D=D_1,D_2,\dots,D_k=D'$ such that either
$\Ext_A^1(D_i,D_{i+1})\ne0$ or $\Ext_A^1(D_{i+1},D_i)\ne0$, for 
$1\le i<k$. More generally, two $A$--modules $M$ and $N$ belong to the
same block if all of their composition factors belong to the same
block.

The natural surjection $\Haff_n\longrightarrow\H_{r,n}(q,\Q)$ shows that if
$D$ and $D'$ are in the same block as $\H_{r,n}(q,\Q)$--modules then they
are in the same block as $\Haff_n$--modules.  The main result of this
paper shows that the blocks of the Ariki--Koike algebras are
determined by the affine Hecke algebra. 

\begin{Theorem}
Suppose that $\F$ is an algebraically closed field and that $q\ne1$.
Let $D$ and $D'$ be irreducible modules for the Ariki--Koike algebra
$\H_{r,n}(q,\Q)$. Then $D$ and $D'$ belong to the same block as
$\H_{r,n}(q,\Q)$--modules if and only if they belong to the same block
as $\Haff_n$--modules.
\end{Theorem}

We also classify the blocks of the Ariki--Koike algebras when $q=1$
and when some of the parameters $Q_1,\dots,Q_r$ are zero. 

By a well--known theorem of Bernstein~\cite[Prop.~3.11]{Lusztig:graded},
the centre of $\Haff_n$ is the set $\F[X_1^\pm,\dots,X_n^\pm]^{\Sym_n}$ of
symmetric Laurent polynomials in $X_1,\dots,X_n$. Consequently, the central
characters of $\Haff_n$ are naturally indexed by $\Sym_n$--orbits of
$(F^\times)^n$. This observation gives a natural combinatorial criterion
for two $\H_{r,n}(q,\Q)$--modules to belong to the same block (see
Theorem~\ref{blocks} for the precise statement), and it is this statement
that we actually prove. We prove Theorem~A by first showing the that blocks
of $\H_{r,n}(q,\Q)$ are the `same' as the blocks of the associated
cyclotomic $q$--Schur algebra. This allows us to use a new characterization
of the blocks in terms of `Jantzen coefficients'
(Proposition~\ref{Jantzen=blocks}).

Observe that the Theorem~A is equivalent to the following property of
the blocks of~$\Haff_n$.

\begin{Corollary}
Suppose that $q\ne1$ and let $D$ and $D'$ be two simple
$\H_{r,n}(q,\Q)$--modules. Then $D$ and~$D'$ belong to the same block
as $\Haff_n$--modules if and only if there exist simple
$\H_{r,n}(q,\Q)$--modules $D=D_1,D_2,\dots,D_k=D'$ such that either
$$\Ext_{\Haff_n}^1(D_i,D_{i+1})\ne0\quad\text{or}\quad
           \Ext_{\Haff_n}^1(D_{i+1},D_i)\ne0,$$
for $1\le i<k$.
\end{Corollary}

In 1999 Grojnowski~\cite{Groj:AKblocks} announced a proof of
Theorem~A.  Using an ingenious argument, what Grojnowski actually
proves is that
$$\Ext_{\Haff_n}^1(D,D')=\Ext_{\H_{r,n}(q,\Q)}^1(D,D')$$
whenever $D\ne D'$ are simple $\H_{r,n}(q,\Q)$--modules.
Unfortunately, as Anton Cox~\cite{Cox:AKblocks} has pointed out, this
is not enough to classify the blocks of the Ariki--Koike algebras. For
example, it could happen that there are no $\Haff_n$--module
extensions between different $\H_{r,n}(q,\Q)$--modules which belong to
the same block as $\Haff_n$--modules. We note that Grojnowski's result
does not follow from Theorem~A.

Lusztig~\cite{Lusztig:graded} introduced a graded, or degenerate,
Hecke algebra for each affine Hecke algebra.
Brundan~\cite{Brundan:degenCentre} has shown that the centre of the
degenerate affine Hecke algebra maps onto the centre of the degenerate
cyclotomic Hecke algebras. This gives a classification of the blocks
of the degenerate cyclotomic and affine Hecke algebras analogous to
our Theorem~A. It should be possible to use the arguments
from this paper to classify the blocks of the degenerate cyclotomic
Hecke algebras of type $G(r,1,n)$ and the associated degenerate
cyclotomic Schur algebras.  All of the combinatorics that we use goes
through without change, however, it is necessary to check that
arguments of \cite{JM:cyc-Schaper} can be adapted to prove a sum
formula for the Jantzen filtrations of the degenerate cyclotomic Schur
algebras. This should be routine (cf.~\cite[\S6]{AMR}), however, we
have not checked the details. 

The outline of this paper is as follows. In the next section we introduce
the Ariki--Koike algebras and the cyclotomic $q$--Schur algebras. Using the
representation theory of these two algebras, we reduce the proof of
Theorem~A to a purely combinatorial problem of showing that two equivalence
relations on the set of multipartitions coincide (Theorem~\ref{blocks}).
The first of these equivalence relations comes from the cyclotomic Jantzen
sum formula~\cite{JM:cyc-Schaper}. The second equivalence relation is
equivalent to the combinatorial criterion which classifies the central
characters the affine Hecke algebras. In section~3 we develop the
combinatorial machinery needed to show that our two equivalence relations
on the set of multipartitions coincide when $q\ne1$ and when the parameters
$Q_1,\dots,Q_r$ are non--zero. Here we are greatly aided by recent work of
Fayers~\cite{Fayers:AKweight,Fayers:multicore} on the `core block' of a
multipartition. Finally, in section~4 we consider the blocks of the
Ariki--Koike algebras with `exceptional' parameters; that is, those
algebras with $q=1$ or with some of the parameters $Q_1,\dots,Q_r$ being
zero. Quite surprisingly, the algebras with exceptional parameters have
only a single block (unless $q=1$ and $r=1$).

\section{Cyclotomic Hecke algebras and Schur algebras}

This section begins by introducing the cyclotomic Hecke algebras and
Schur algebras. We then reduce the proof of Theorem~A to a
purely combinatorial statement which amounts to showing that two
equivalence relations on the set of multipartitions coincide.

\subsection{Ariki--Koike algebras}
Let $\F$ be a field of characteristic
$p\in\{2,3,\dots\}\cup\{\infty\}$ and fix positive integers~$n$
and~$r$. Suppose that $q,Q_1,\ldots Q_r$ are elements of $\F$ such
that $q$ is invertible and let $\Q=(Q_1,\dots,Q_r)$. The
\emph{Ariki--Koike algebra} $\h_{r,n} =\h_{r,n}(q,\Q)$ is the unital
associative $\F$--algebra with generators $T_0,T_1,\ldots,T_{n-1}$ and
relations
\begin{align*}
(T_i+q)(T_i-1) & = 0, & 1 & \leq i \leq n-1, \\
(T_0-Q_1)\ldots(T_0-Q_r) & =0, && \\
T_i T_j & = T_j T_i, & 0 &\leq i <j-1 \leq n-2, \\
T_i T_{i+1} T_i & = T_{i+1} T_i T_{i+1}, & 1 & \leq i \leq n-2, \\ 
T_0 T_1 T_0 T_1 & = T_1 T_0 T_1 T_0. &&
\end{align*}
Define $e \geq 2$ to be minimal such that $1+q+\ldots+ q^{e-1}=0 \in
\F$. Then $e \in \{2,3,\ldots\} \cup \{\infty\}$. 
Note that $e=p$ if and only if $q=1$. If $e\ne p$ and $p$ is finite
then $p$ does not divide~$e$.

Recall that a partition $\lambda=(\lambda_1,\lambda_2,\dots)$ of~$n$
is a weakly decreasing sequence of non--negative integers which sum to
$|\lambda|=n$. An $r$--multipartition of $n$, or more simply a
multipartition, is an ordered $r$--tuple
$\bl=\(\lambda^{(1)},\dots,\lambda^{(r)})$ of partitions with
$|\bl|=|\la^{(1)}|+\dots+|\la^{(r)}|=n$. Let $\Lamb$ be the set of
multipartitions of~$n$.  We regard a partition as a multipartition
with one component, so any subsequent definition concerning
multipartitions specializes to a corresponding definition for
partitions.

The set of multipartitions is naturally
ordered by \emph{dominance} where $\bl\gedom\bmu$ if
$$\sum_{t=1}^{s-1}|\lambda^{(t)}|+\sum_{j=1}^i\lambda^{(s)}_j
   \ge\sum_{t=1}^{s-1}|\mu^{(t)}|+\sum_{j=1}^i\mu^{(s)}_j$$
for $s=1,2,\dots,r$ and all $i\ge1$. We write $\bl\gdom\bmu$ if
$\bl\gedom\bmu$ and $\bl\ne\bmu$.

The Ariki--Koike algebra $\H_{r,n}$ is a cellular
algebra~\cite{GL,DJM:cyc}. The cell modules of $\H_{r,n}$ are
indexed by the multipartitions of $n$. The cell module indexed by the
multipartition $\bl$ is the \emph{Specht module}~$S(\bl)$. By the
theory of cellular algebras~\cite{GL,M:ULect}, there is an
$\H_{r,n}$--invariant bilinear form $\<\, , \, \>_\bl$ on the Specht 
module $S(\bl)$, so the radical $\rad S(\bl)=\set{x\in
S(\bl)|\<x,y\>_\bl=0\text{ for all }y\in S(\bl)}$ is an
$\H_{r,n}$--submodule of~$S(\bl)$. Set $D(\bl)=S(\bl)/\rad S(\bl)$. Then the
non--zero $D(\bl)$ give a complete set of pairwise non--isomorphic simple
$\H_{r,n}$--modules.

The theory of cellular algebras gives us the following fact which
is vital for this paper because it allows us work with Specht modules
rather than with the simple $\H_{r,n}$--modules.

\begin{lemma}[\protect{Graham--Lehrer~\cite[3.9.8]{GL},
    \cite[Cor.~2.2]{M:ULect}}]
\label{cell reduction}
Suppose that $\bl$ is a multipartition. Then all of the composition
factors of~$S(\bl)$ belong to the same block.
\end{lemma}

Thus we can talk of  the block of $\H_{r,n}$ which contains the Specht
module $S(\bl)$. 

\subsection{Cyclotomic $q$--Schur algebras}
Rather than working with Specht modules to classify the blocks we want
to work with Weyl modules. To this end let 
$\set{L_1^{a_1}\dots L_n^{a_n}T_w|0\le a_i<r\text{ and }w\in\Sym_n}$
be the Ariki-Koike basis of $\H_{r,n}$~\cite[Prop.~3.4]{AK}. That is,
$L_1=T_0$ and $L_{i+1}=q^{1-i}T_iL_iT_i$, for $1\leq i<n$, and if
$w\in\Sym_n$ then $T_w=T_{i_1}\dots T_{i_k}$ whenever
$w=(i_1,i_1+1)\dots(i_k,i_k+1)$ with $k$ minimal (so this is a reduced
expression of $w$). For each multipartition $\bl$ define
$$m_\bl=\prod_{s=1}^r\prod_{k=1}^{|\bl^{(1)}|+\dots+|\bl^{(s-1)}|}
             (L_k-Q_s)\cdot \sum_{w\in\Sym_\bl}T_w,$$
where $\Sym_\bl=\Sym_{\bl^{(1)}}\times\dots\times\Sym_{\bl^{(r)}}$
is the Young subgroup of $\Sym_n$ associated to~$\bl$. The
\emph{cyclotomic $q$--Schur algebra} is the endomorphism algebra
$$\Sch_{r,n}=\Sch_{r,n}(q,\Q)
  =\End_{\H_{r,n}}\Big(\bigoplus_{\bl\in\Lamb}m_\bl\H_{r,n}\Big).$$
We remark that this variant of the cyclotomic $q$--Schur algebra is Morita
equivalent to one of the algebras introduced in~\cite{DJM:cyc}.
The representation theory of $\Sch_{r,n}$ is discussed
in~\cite{M:cyclosurv}.

The cyclotomic $q$--Schur algebra $\Sch_{r,n}$ is a quasi--hereditary
cellular algebra.  The cell modules of $\Sch_{r,n}$ are the \emph{Weyl
modules} $\Delta(\bl)$, for $\bl\in\Lamb$. For each
$\bl\in\Lamb$, there is a non--zero simple module
$L(\bl)=\Delta(\bl)/\rad\Delta(\bl)$. Just as with 
Lemma~\ref{cell reduction}, the theory of cellular algebras tells us
the following.

\begin{lemma}[\protect{Graham--Lehrer~\cite[3.9.8]{GL},
    \cite[Cor.~2.2]{M:ULect}}]
\label{W cell reduction}
Suppose that $\bl$ is a multipartition. Then all of the composition
factors of~$\Delta(\bl)$ belong to the same block.
\end{lemma}

The next result shows that in order to classify the blocks of
$\H_{r,n}$ it is enough to consider the blocks of~$\Sch_{r,n}$. In
fact, this is an easy consequence of double centralizer theory.

Let $A$ be a finite dimensional algebra over a field. Then $A$ decomposes
in a unique way as a direct sum of indecomposable two-sided ideals
$\H=B_1\oplus\dots\oplus B_d$. Recall that two simple $A$-modules~$D$ and
$D'$ are in the same block if there exist simple modules
$D_1=D,D_2,\dots,D_k=D'$ such that $\Ext^1_A(D_i,D_{i+1})\ne0$ or
$\Ext^1_A(D_{i+1},D_i)\ne0$, for $1\le i<k$. As $\Ext^1_A$ classifies
non-trivial extensions, it follows that two simple modules $D$ and $D'$
belong to the same block if and only if $D$ and $D'$ are both composition
factors of $B_j$ or, equivalently, that $D=DB_j$ and $D'=D'B_j$, for
some~$j$. Abusing terminology, we call the indecomposable subalgebras
$B_1,\dots,B_d$ the blocks of~$A$ and we say that an $A$-module $M$ belongs
to the block $B_j$ if $MB_j=M$.  Using an idempotent argument
(cf.~\cite[Theorem~56.12]{C&R:2} it is now easy to show that two
indecomposable $A$--modules $P$ and $Q$ belong to the same block if and
only if they are in the same \emph{linkage class}; that is, there exist
indecomposable modules $P_1=P,\dots,P_l=Q$ such that $P_i$ and~$P_{i+1}$
have a common irreducible composition factor, for $i=1,\dots,l-1$. 

By Lemma~\ref{cell reduction} and the last paragraph, that two Specht
modules $S(\bl)$ and $S(\bmu)$ belong to the same block if and only if
there exist multipartitions $\bl_1=\bl,\dots,\bl_k=\bmu$ such that
$S(\bl_i)$ and $S(\bl_{i+1})$ have a common composition factor, for $1\le
i<k$. Similarly, two Weyl modules belong to the same block if and only if
they are in the same linkage class. We will use this characterization of
the blocks of $\H$ and~$\Sch_{r,n}$ below without mention.

\begin{proposition}\label{same}
Let $\bl$ and $\bmu$ be multipartitions of $n$. Then $S(\bl)$
and $S(\bmu)$ are in the same block as $\h_{r,n}$--modules if and only
if $\Delta(\bl)$ and $\Delta(\bmu)$ are in the same block as
$\Sch_{r,n}$--modules.
\end{proposition}

\begin{proof} Suppose first that $S(\bl)$ and $S(\bmu)$ are in the
same block. By Lemma~\ref{cell reduction} all of the composition
factors of $S(\bl)$ belong to the same block. Therefore, by the
remarks above, it is enough to consider the case when
$D(\bmu)\ne0$ and $D(\bmu)$ is a composition factor of $S(\bl)$. By a
standard Schur functor argument~\cite[Prop.~2.17]{JM:cyc-Schaper},
$[\Delta(\bl){:}L(\bmu)]=[S(\bl){:}D(\bmu)]\ne0$. Therefore,
$\Delta(\bl)$ and $\Delta(\bmu)$ are in the same block.  Note that
this implies that $\Sch_{r,n}$ cannot have more blocks (that is,
indecomposable subalgebras) than $\H_{r,n}$.

To prove the converse let
$M=\bigoplus_{\bl\in\Lamb}m_\bl\H_{r,n}$ and suppose that
$\H_{r,n}=B_1\oplus\dots\oplus B_k$ is the unique decomposition of
$\H_{r,n}$ into a direct sum of indecomposable subalgebras. Then 
$$M=M\H_{r,n}=MB_1+\dots+ MB_k.$$
In fact, this sum is direct because, by definition, 
$MB_i\cap MB_j=\emptyset$ if $i\ne j$, and  $MB_i\ne0$ since $\H_{r,n}$
is a submodule of $M$. Therefore,
\begin{align*}
\Sch_{r,n}&=\End_{\H_{r,n}}(M)
           =\End_{\H_{r,n}}\(MB_1\oplus\dots\oplus MB_k\)\\
	  &=\bigoplus_{1\le i,j\le k}\Hom_{\H_{r,n}}(MB_i,MB_j)
           =\bigoplus_{i=1}^k\End_{\H_{r,n}}(MB_i),
\end{align*}
where the last equality follows because $B_i$ and $B_j$ have no
common irreducible constituents if $i\ne j$. Consequently, $\Sch_{r,n}$
has at least as many blocks as $\H_{r,n}$.

Combining the last two paragraphs proves the proposition.
\end{proof}

Thus, to prove Theorem~A it suffices to determine when two Weyl
modules are in the same block. The advantage of working with Weyl
modules is shown in Lemma~\ref{vanishing} below. Before we can state
this result we need some notation. 

If $A$ is an algebra let $K_0(A)$ be the Grothendieck group of finite
dimensional $A$--modules and if $M$ is a $A$--module let $[M]$ be its
image in $K_0(A)$. In particular, the Grothendieck group
$K_0(\Sch_{r,n})$ of $\Sch_{r,n}$ is the free $\Z$--module with basis
$\set{[L(\bl)]|\bl\in\Lamb}$.  The images
$\set{[\Delta(\bl)]|\bl\in\Lamb}$ of the Weyl modules give a
second basis of $K_0(\Sch_{r,n})$ since $[\Delta(\bl){:}L(\bl)]=1$ and $[\Delta(\bl){:}L(\bmu)]>0$ only if $\bl \gedom \bmu$, for all
$\bl,\bmu\in\Lamb$ (see~\cite{GL}). Hence, we have the following.

\begin{lemma}\label{vanishing}
    Suppose that $a_\bl\in\Z$. Then
    $\sum_\bl a_\bl[\Delta(\bl)]=0$ in $K_0(\Sch_{r,n})$ if and only if
    $a_\bl=0$ for all~$\bl\in\Lamb$.
\end{lemma}

Note that, in general, there can exist non--zero integers $a_\bl\in\Z$
such that $\sum_\bl a_\bl[S(\bl)]=0$. This follows because 
$K_0(\H_{r,n})$ is a free $\Z$--module of rank
$L=\#\set{\bl\in\Lamb|D(\bl)\ne0}$ and $L=\#\Lamb$ (if and) only if
$\H_{r,n}$ is semisimple.

\subsection{The cyclotomic Jantzen sum formula}
The next step is to recall (a special case of) the machinery of the
cyclotomic Jantzen sum formula~\cite{JM:cyc-Schaper}. Let $t$ be an
indeterminate over $\F$ and let $\O=\F[t,t^{-1}]_\pi$ be the
localization of $\F[t,t^{-1}]$ at the prime ideal $\pi=\<t-1\>$. Let
$\Sch_\O=\Sch_\O(qt,\X)$ be the cyclotomic Schur algebra over~$\O$ with
parameters~$qt$ and $\X=(X_1,\dots,X_r)$ where
$$X_a=\begin{cases}
    Q_at^{na},&\text{if }Q_a\ne0,\\
    (t-1)t^{na},&\text{if }Q_a=0.
\end{cases}$$ 
Consider $\F$ as an $\O$--module by letting $t$ act on~$\F$ as multiplication
by~$1$. Then $\Sch_{r,n}\cong\Sch_\O\otimes_\O\F$, since $\Sch_\O$ is
free as an $\O$--module by~\cite[Theorem~6.6]{DJM:cyc}. The algebra
$\Sch_\O\otimes_\O\F(t)$ is split semisimple by Schur--Weyl
duality~\cite[Theorem~5.3]{M:cyclosurv} and Ariki's criterion for the
semisimplicity for~$\H_{r,n}$~\cite{Ariki:ss}. Thus we are in the
general setting considered in~\cite[\S4]{JM:cyc-Schaper}.

Let $\nu_\pi$ be the $\pi$--adic evaluation map on $\O^\times$; thus,
$\nu_\pi(f(t))=k$ if $k\ge0$ is maximal such that $(t-1)^k$ divides
$f(t)\in\F[t,t^{-1}]$. Let $\Delta_\O(\bl)$ be the Weyl module of
$\Sch_\O$ indexed by the multipartition $\bl\in\Lamb$. Recall that
$\Delta_\O(\bl)$ carries a bilinear form $\<\ ,\ \>_\bl$ by the
general theory of cellular algebras. For each integer $i\ge0$ 
define
$$\Delta_\O(\bl)_i=\set{x\in\Delta_\O(\bl)|\nu_\pi(\<x,y\>)\ge i
                      \text{ for all }y\in\Delta_\O(\bl)}.$$
Finally, let 
$\Delta(\bl)_i=\(\Delta_\O(\bl)_i+\pi\Delta_\O(\bl)\)/\pi\Delta_\O(\bl)$. Then
$$\Delta(\bl)=\Delta(\bl)_0\supset\Delta(\bl)_1
             \supseteq\Delta(\bl)_2\supseteq\dots$$
is a \emph{Jantzen filtration} of the $\Sch_{r,n}$--module $\Delta(\bl)$.
Then $\Delta(\bl)_k=0$ for $k\gg0$ since $\Delta(\bl)$ is finite
dimensional.

To describe the Jantzen filtration of $\Delta(\bl)$ we need some
combinatorics. The \emph{diagram} of a multipartition $\bl$ is the set
$[\bl]=\set{(i,j,a)|1\le j\le\lambda^{(a)}_i\text{ and }1\le a\le r}$.  A
\emph{node} is any ordered triple $(i,j,a)$ 
in~$\N\times\N\times\{1,\dots,r\}$.  For example, all of the elements of
$[\bl]$ are nodes.  

Each node $x=(i,j,a)\in[\bl]$ determines a \emph{rim hook} 
$$r^\bl_x=\set{(k,l,a)\in[\bl]|k\ge i, l\ge j
                     \text{ and }(k+1,l+1,a)\notin[\bl]}.$$
We say that $r^\bl_x$ is an \emph{$h$--rim hook} if $h=|r^\bl_x|$.  Let
$i'$ be maximal such that $(i',j,a)\in[\bl]$; so $i'$ is the length of
column~$j$ of $\lambda^{(a)}$. Then $f^\bl_x=(i',j,a)\in[\bl]$ is the
\emph{foot} of $r^\bl_x$ and $r^\bl_x$ has \emph{leg length}
$\leg(r^\bl_x)=i'-i$. Similarly, the node $(i,\lambda^{(a)}_j,a)$ is
the \emph{hand} of $r^\bl_x$.  If $x\in[\bl]$ let
$\bl{\setminus}r^\bl_x$ be the multipartition with diagram
$[\bl]{\setminus} r^\bl_x$. We say that $\bl{\setminus}r^\bl_x$ is the
multipartition obtained by \emph{unwrapping} the rim hook~$r^\bl_x$
from~$\bl$, and that $\bl$ is the multipartition obtained from
$\bl{\setminus}r^\bl_x$ by \emph{wrapping} on the rim hook~$r^\bl_x$.

Define the $\O$--residue of the node $x=(i,j,a)$ to be
$\res_\O(x)=(qt)^{j-i}X_a$.

\begin{definition}\label{Jlambdamu}
Suppose that $\bl=(\lambda^{(1)},\dots,\lambda^{(r)})$ and
$\bmu=(\mu^{(1)},\dots,\mu^{(r)})$ are multipartitions of~$n$. 
The \emph{Jantzen coefficient} $J_{\bl\bmu}$ is the integer
$$J_{\bl\bmu}=\begin{cases}
   \quad\Sum_{x\in[\bl]}
      \sum_{\substack{y\in[\bmu]\\\relax
             [\bmu]{\setminus}r^\bmu_y=[\bl]{\setminus}r^\bl_x}}
              (-1)^{\leg(r^\bl_x)+\leg(r^\bmu_y)}
	      \nu_\pi\(\res_\O(f^\bl_x)-\res_\O(f^\bmu_y)\),
   &\text{if }\bl\gdom\bmu,\\
  \quad0,&\text{otherwise.}
\end{cases}
$$
\end{definition}

The Jantzen coefficient $J_{\bl\bmu}$ depends on the choices of $\F$,
$q$ and $\Q=(Q_1,\dots,Q_r)$. In fact, we will see that $J_{\bl\bmu}$ depends only on
$p$,~$e$ and $\Q$.  By definition $J_{\bl\bmu}$ is an integer which is
determined by the combinatorics of multipartitions.  The definition of
$J_{\bl\bmu}$ is reasonably involved, however, it turns out that these
integers are computable. In sections~3 and~4 we give simpler formulae for
the Jantzen coefficients.

\begin{theorem}[James and Mathas~\cite{JM:cyc-Schaper}, Theorem~4.3]
\label{Jantzen}
Suppose that $\bl$ is a multipartition of $n$. Then
$$\sum_{i>0}\,[\Delta(\bl)_i]
  =\sum_{\bmu\in\Lamb} J_{\bl\bmu}\, [\Delta(\bmu)]$$
in $K_0(\Sch_{r,n})$.
\end{theorem}

For  multipartitions $\bl$ and  $\bmu$ in $\Lamb$ let
$d_{\bl\bmu}=[\Delta(\bl){:}L(\bmu)]$ be the number of composition
factors of~$\Delta(\bl)$ which are isomorphic to $L(\bmu)$. Define
$$J_{\bl\bmu}'=\sum_{\substack{\bnu\in\Lamb\\\bl\gdom\bnu\gedom\bmu}}
          J_{\bl\bnu}d_{\bnu\bmu}.$$
By Theorem~\ref{Jantzen}, $J_{\bl\bmu}'$ is the composition
multiplicity of the simple module $L(\bmu)$ in
$\bigoplus_{i>0}\Delta(\bl)_i$. Therefore, $J_{\bl\bmu}'\ge0$, for all
$\bl,\bmu\in\Lamb$. As $\Delta(\bl)_1=\rad\Delta(\bl)$ we obtain the following.

\begin{corollary}\label{adjacent}
Suppose that $\bl\ne\bmu$ are multipartitions of $n$. Then
$d_{\bl\bmu}\le J_{\bl\bmu}'$ and, moreover, $d_{\bl\bmu}\ne0$ if and
only if~$J_{\bl\bmu}'\ne0$.
\end{corollary}

We now use Theorem~\ref{Jantzen} to classify the blocks of $\Sch_{r,n}$.

\begin{definition}
Suppose that $\bl,\bmu\in\Lamb$. Then $\bl$ and $\bmu$ are
\emph{Jantzen equivalent}, and we write $\bl\sim_J\bmu$, if
there exists a sequence of multipartitions $\bl_0=\bl,\bl_1,\dots,\bl_k=\bmu$ 
such that either
$$J_{\bl_i\bl_{i+1}}\ne0\qquad\text{or}\qquad J_{\bl_{i+1}\bl_i}\ne0,$$
for $0\le i < k$.
\end{definition}

Jantzen equivalence gives us our first combinatorial characterization
of the blocks of $\Sch_{r,n}$.

\begin{proposition}\label{Jantzen=blocks}
Suppose that $\bl,\bmu\in\Lamb$. Then $\Delta(\bl)$ and
$\Delta(\bmu)$ belong to the same block as $\Sch_{r,n}$--modules if
and only if~$\bl\sim_J\bmu$.
\end{proposition}

\begin{proof}
We first show that $\Delta(\bl)$ and $\Delta(\bmu)$ belong to the same
block whenever $\bl\sim_J\bmu$. By definition $\Delta(\bl)_i$ is a
submodule of $\Delta(\bl)$ for all~$i$, so  all of the composition
factors of $\sum_{i>0}\Delta(\bl)_i$ belong to the same block as $\Delta(\bl)$ by Lemma~\ref{W cell reduction}. Consequently, all of the composition
factors of the virtual module $\sum_\bnu J_{\bl\bnu}[\Delta(\bnu)]$
belong to the same block. Let $\Lambda'$ be the set of
multipartitions $\bnu$ such that $\Delta(\bnu)$ is not in the same
block as $\Delta(\bl)$. Then we have
$\sum_{\bnu\in\Lambda'}J_{\bl\bnu}[\Delta(\bnu)]=0$. Hence,
$J_{\bl\bnu}=0$ whenever $\bnu\in\Lambda'$ by
Lemma~\ref{vanishing}. It follows that $\Delta(\bl)$ and
$\Delta(\bmu)$ belong to the same block whenever $\bl\sim_J\bmu$.

To prove the converse it is sufficient to show that $\bl\sim_J\bmu$
whenever $d_{\bl\bmu}\ne0$. Hence, by Corollary~\ref{adjacent} we must
show that $\bl\sim_J\bmu$ whenever $J_{\bl\bmu}'\ne0$. However, if
$J_{\bl\bmu}'\ne0$ then we can find a multipartition~$\bnu_1$ such
that $J_{\bl\bnu_1}\ne0$, $d_{\bnu_1\bmu}\ne0$ and
$\bl\gdom\bnu_1\gedom\bmu$. Consequently, $\bl\sim_J\bnu_1$.  If
$\bnu_1\ne\bmu$ then $J_{\bnu_1\bmu}'\ne0$ by Corollary~\ref{adjacent}
since $d_{\bnu_1\bmu}\ne0$. Therefore, we can find a multipartition
$\bnu_2$ such that $J_{\bnu_1\bnu_2}\ne0$, $d_{\bnu_2\bmu}\ne0$ and
$\bnu_1\gdom\bnu_2\gedom\bmu$. Continuing in this way we can find
multipartitions $\bnu_0=\bl,\bnu_1,\dots,\bnu_k=\bmu$ such that
$J_{\bnu_{i-1}\bnu_i}\ne0$, $d_{\bnu_i\bmu}\ne0$, for $0< i<k$, 
and $\bl\gdom\bnu_1\gdom\dots\gdom\bnu_k=\bmu$. Note that
we must have $\bnu_k=\bmu$ for some $k$ since $\Lamb$ is finite. Therefore,
$\bl\sim_J\bnu_1\sim_J\dots\sim_J\bnu_k=\bmu$ as required.
\end{proof}

\begin{remark}
Under some very mild technical assumptions (see, for example,
\cite[\S4.1]{McNinch:Howe}), Jantzen filtrations can be defined for
the standard modules of an arbitrary quasi--hereditary algebra.  The
argument of Proposition~\ref{Jantzen=blocks} is completely generic: it
shows that the blocks of a quasi--hereditary algebra are determined by
the `Jantzen coefficients'.
\end{remark}

\begin{remark}Without using the cyclotomic $q$--Schur algebras it is
    not clear that Jantzen equivalence determines the blocks of
    $\H_{r,n}$. Applying the Schur functor to Theorem~\ref{Jantzen}
    gives an analogous description of the Jantzen filtration of the
    Specht modules: $\sum_{i>0}[S(\bl)_i]=\sum_\bmu
    J_{\bl\bmu}[S(\bmu)]$.  The problem is that, \textit{a priori}, 
    the composition factors of 
    $\bigoplus_\bmu J_{\bl\bmu}S(\bmu)$ could belong to different
    blocks because the analogue of Lemma~\ref{vanishing} fails for
    Specht modules.
\end{remark}

\subsection{A second combinatorial characterization of the blocks}
Proposition~\ref{Jantzen=blocks} completely determines the blocks of
$\Sch_{r,n}$, and hence the blocks of $\H_{r,n}$. Unfortunately, it is
not obvious when two multipartitions are Jantzen equivalent. 

The \emph{residue} of the node $x=(i,j,a)$ is 
$$\res(x)=\begin{cases}
    q^{j-i}Q_a,&\text{if $q\ne1$ and $Q_a\ne0$},\\
    (\bar{j-i},Q_a),&\text{if $q=1$ and $Q_a\ne Q_b$ for $b\ne a$},\\
    Q_a,&\text{otherwise,}
\end{cases}$$
where $\bar z=z\pmod p$ for $z\in\Z$ (if $p=\infty$ we set $\bar
z=z$). Let 
$$\Res=\set{\res(x)|x\in[\bl]\text{ for some }\bl\in\Lamb}$$ 
be the set of all possible residues.  For any multipartition
$\bl\in\Lamb$ and $f\in\Res$ define
$$C_f(\bl)=\#\set{x\in[\bl]|\res(x)=f}.$$  We can now define our
second combinatorial equivalence relation on~$\Lamb$.

\begin{definition}Suppose that $\bl$ and $\bmu$ are multipartitions.
    Then $\bl$ and $\bmu$ are \emph{residue equivalent}, and we write
    $\bl\sim_C\bmu$, if $C_f(\bl)=C_f(\bmu)$ for all $f\in\Res$.
\end{definition}

It is easy to determine if two multipartitions are residue equivalent,
so the next result gives an effective characterization of the blocks
of the algebras $\H_{r,n}$ and $\Sch_{r,n}$.

{\samepage
\begin{theorem}\label{blocks}
Suppose that $\bl$ and $\bmu$ are multipartitions of~$n$. Then the 
following are equivalent.
\begin{enumerate}
\item $S(\bl)$ and $S(\bmu)$ belong to the same block as 
    $\H_n(\Q)$--modules.
\item $\Delta(\bl)$ and $\Delta(\bmu)$ belong to the same block as
    $\Sch_{r,n}(\Q)$--modules.
\item $\bl\sim_J\bmu$.
\item $\bl\sim_C\bmu$.
\end{enumerate}
\end{theorem}
}

By Proposition~\ref{same} and Proposition~\ref{Jantzen=blocks}, (a),
(b) and (c) are equivalent. Therefore, to prove the theorem it is
enough to prove that $\bl\sim_J\bmu$ if and only if $\bl\sim_C\bmu$.
The proof of this fact is given in sections~3 and~4. It turns out
that, combinatorially, these equivalence relations depend very much on
whether or not~$q=1$ and whether or not some of the parameters
$Q_1,\dots,Q_r$ are zero. The following result allows us to treat
these cases separately.

\begin{theorem}%
[Dipper and Mathas~\cite{DM:Morita}, Theorem~1.5 and Corollary~5.7]%
\leavevmode\newline\label{morita}%
Suppose that $\Q=\Q_1\coprod\Q_2\coprod\dots\coprod\Q_\kappa$ is a
partition of $\Q$ such that $q^cQ_a\in\Q_\alpha$ only if
$Q_a\in\Q_\alpha$, for~$c\in\Z$, $1\le a\le r$ and
$1\le\alpha\le\kappa$. Set $r_\alpha=|\Q_\alpha|$, for 
$1\le \alpha\le\kappa$. Then $\Sch_{r,n}(\Q)$ is Morita equivalent~to
the algebra
$$\bigoplus_{\substack{n_1,\dots,n_\kappa\ge0\\n_1+\dots+n_\kappa=n}}
     \Sch_{r_1,n_1}(\Q_1)\boxtimes\Sch_{r_2,n_2}(\Q_2)\boxtimes
        \dots\boxtimes \Sch_{r_\kappa,n_\kappa}(\Q_\kappa).$$ 
Moreover, if $\Q_\alpha=\{Q_{i^\alpha_1},\dots,Q_{i^\alpha_{r_\alpha}}\}$, for
$1\le\alpha\le\kappa$, then the Morita equivalence is induced by the map
$\Delta(\bl)\mapsto\Delta(\bl_1)\boxtimes\dots\boxtimes\Delta(\bl_\kappa)$,
where $\bl_\alpha=(\lambda^{(i^\alpha_1)},\dots,\lambda^{(i^\alpha_{r_\alpha})})$,
for $1\le\alpha\leq\kappa$ and $\bl\in\Lamb$.
\end{theorem}

There is an analogous result for the Ariki--Koike algebra $\H_{r,n}$; see
\cite[Theorem~1.1]{DM:Morita}.

Theorem~\ref{morita} says that the blocks of $\H_{r,n}(\Q)$ and
$\Sch_{r,n}(\Q)$ depend only on the $q$--orbits of the parameters;
that is, the orbits of the parameters under multiplication by~$q$.
Further, by Theorem~\ref{morita} it is enough to consider the case
where $\Q$ is contained in a single $q$--orbit to prove
Theorem~\ref{blocks}. Hence, by rescaling $T_0$, if necessary, we can
assume that the parameters $Q_1,\dots,Q_r$ are all are all zero or
that they are all powers of~$q$. More explicitly, we can assume that
either $Q_a=0$, or that there exist integers
$c_1,\dots,c_r$ such that $Q_a=q^{c_a}$, for $1\le a\le r$.
Consequently, to prove Theorem~\ref{blocks} we are reduced to
considering the following mutually exclusive cases:
\refstepcounter{equation}\label{cases}
$$\begin{tabular}[t]{llp{60mm}}
    Case 1.& $q\ne1$ and $Q_a=q^{c_a}$, for $1\le a\le r$.\\
    Case 2.& $r=1$ and $q=1$.&\\
    Case 3.& $r>1$, $q=1$ and $Q_1=\dots=Q_r=1$.\\
    Case 4.& $r>1$, $q=1$ and $Q_1=\dots=Q_r=0$.\\
    Case 5.& $r>1$, $q\ne1$ and $Q_1=\dots=Q_r=0$.
\end{tabular}\leqno(\theequation)$$
Note that $\H=\H_{1,n}$ is independent of $Q_1$ when $r=1$.

The proof of Theorem~\ref{blocks} for case~1 is given in section~3.
Cases~2--5 are considered in section~4 using similar, but easier,
arguments. Given a node $x=(i,j,a)$ note that $\res(x)=q^{j-i}Q_a$ in
case~1, $\res(x)=(\bar{j-i},Q_1)$ in case~2 and $\res(x)=Q_a$ in the
other three cases.

The basic strategy for proving Theorem~\ref{blocks} for each of these
five cases is the same, however, the proof breaks up into three cases
because the combinatorics of residue equivalence is different for
case~1, case~2 and cases~3--5. Fayers has pointed out that the
Ariki--Koike algebras in cases~3 and~4 are isomorphic via the algebra
homomorphism determined by $T_0\mapsto (T_0-1)$ and $T_i\mapsto T_i$,
for $1\le i<n$, so we do not actually need to consider case~4.

\subsection{The blocks of the affine Hecke algebra}Assuming
Theorem~\ref{blocks} we now prove Theorem~A from the
introduction. 

As the centre $Z(\Haff_n)$ of $\Haff_n$ is the set of symmetric
Laurent polynomials in $X_1,\dots,X_n$, the central characters of
$\Haff_n$ are indexed by $\Sym_n$--orbits of $(\F^\times)^n$. More
precisely, if $\gamma\in(\F^\times)^n/\Sym_n$ then the central
character $\chi_\gamma$ is given by evaluation at $\gamma$. 

By Lemma~\ref{cell reduction}, all of the composition factors of the
Specht module $S(\bl)$ belong to the same block as
$\H_{r,n}$--modules. Therefore, all of the composition factors of
$S(\bl)$ belong to the same block as $\Haff_n$--modules.  We need to
know the central characters of the Specht modules. 

\begin{lemma}\label{cchars}
Suppose that $q\ne1$ and that $D(\bl)\ne0$, for some multipartition
$\bl\in\Lamb$. Then $f(X)\in\Z(\Haff_n)$ acts on $D(\bl)$ as
multiplication by $f(\gamma)$, where
$\gamma=\(\!\res(x_1),\res(x_2),\dots,\res(x_n)\)$ and
$[\bl]=\{x_1,\dots,x_n\}$ (in any order).
\end{lemma}

\begin{proof}
As all of the composition factors of $S(\bl)$ belong to the same block
as $D(\bl)$, $f(X)$ acts on $S(\bl)$ and on $D(\bl)$ as multiplication
by the same scalar. By \cite[Prop.~3.7]{JM:cyc-Schaper} this scalar is
given by evaluating the polynomial $f(X)$ at
$\(\!\res(x_1),\res(x_2),\dots,\res(x_n)\)$.
\end{proof}

\begin{theorem}\label{B'}
Suppose that $q\ne1$ and that $\F$ is algebraically closed. Then two
simple $\Haff_n$--modules~$D$ and~$D'$ belong to the same block if and
only if they have the same central character.
\end{theorem}

\begin{proof}
Any two simple modules in the same block have the same central
character. Conversely, suppose that $D$ and $D'$ are simple
$\Haff_n$--modules which have the same central character. Let
$(X_1-Q_1)\dots(X_1-Q_s)$ and $(X_1-Q_{s+1})\dots(X_1-Q_r)$,
respectively, be the minimal polynomials for $X_1$ acting on $D$ and
$D'$. (Note that $Q_1,\dots,Q_r$ are non--zero since $X_1,\dots,X_n$
are invertible.) Then $D$ and $D'$ are both simple modules for the
Ariki--Koike algebra $\H_{r,n}$ with parameters $Q_1,\dots,Q_r$.
Therefore, $D\cong D(\bl)$ and $D'\cong D(\bmu)$ for some
multipartitions $\bl,\bmu\in\Lamb$.  By assumption, $D$ and $D'$ have
the same central characters. The central character of $D(\bl)$ is
uniquely determined by the multiset of residues
$\set{\res(x)|x\in[\bl]}$ by Lemma~\ref{cchars}. Similarly, the
central character of $D(\bmu)$ is determined by the multiset
$\set{\res(x)|x\in[\bmu]}$. Hence, $C_f(\bl)=C_f(\bmu)$, for all
$f\in\Res$. Therefore, $\bl\sim_C\bmu$, so $D\cong D(\bl)$ and
$D'=D(\bmu)$ are in the same block as $\H_{r,n}$--modules by
Theorem~\ref{blocks}. Hence,~$D$ and $D'$ are in the same block as
$\Haff_n$--modules.  
\end{proof}

Theorem~\ref{B'} is not new. We are grateful to Iain Gordon for
pointing out that the classification of the blocks of $\Haff_n$ by
central characters is an immediate corollary of a general result of
M\"uller~\cite[Theorem~7]{Mueller} since $\Haff_n$ is finite
dimensional over its centre. See also \cite[III.9]{BrownGoodearl}.

Combining Theorem~\ref{blocks} and Theorem~\ref{B'} we obtain a
more descriptive version of Theorem~A.

\begin{corollary}[Theorem~A]\label{combinatorial}
Suppose that $\F$ is an algebraically closed field, $q\ne1$ and
that the parameters $Q_1,\dots,Q_r$ are non--zero.  Let $\bl$ and
$\bmu$ be multipartitions in $\Lamb$ with $D(\bl)\ne0$ 
and~$D(\bmu)\ne0$. Then the following are equivalent:
\begin{enumerate}
\item $D(\bl)$ and $D(\bmu)$ belong to the same block as
$\H_{r,n}$--modules.
\item $D(\bl)$ and $D(\bmu)$ belong to the same block as
$\Haff_n$--modules.
\item $D(\bl)$ and $D(\bmu)$ have the same central character as
$\Haff_n$--modules.
\item $\bl\sim_C\bmu$.
\end{enumerate}
\end{corollary}

\section{Combinatorics}
In this section, we prove $\bl\sim_J\bmu$ if and only if
$\bl\sim_C\bmu$, for $\bl,\bmu\in\Lamb$ in the cases when $q\ne1$ and
all of the parameters $Q_1,\dots,Q_r$ are powers of~$q$. This is
case~1 of~(\ref{cases}). The basic idea is to reduce the comparison of
the Jantzen and residue equivalence relations to the case where the
multipartitions $\bl$ and $\bmu$ are both `cores'. The complication is
that, unlike for partitions (the case $r=1$), we do not have a good
notion of `core' for multipartitions when $r>1$. We circumvent this
difficulty using ideas of
Fayers~\cite{Fayers:AKweight,Fayers:multicore}.

As we are assuming that the parameters $Q_1,\dots,Q_r$ are all powers
of~$q$, there exist integers $c_1,\dots,c_r$ such that $Q_a=q^{c_a}$,
for $1\le a\le r$.  The sequence $\c=(c_1,\dots,c_r)$ is called the
\emph{multi--charge} of $\Q$. 

Now that $\Q$ is contained in a single $q$--orbit, we redefine
the \emph{residue} of a node $x=(i,j,a)$ to be
$$\res(x)=(j-i+c_a)\pmod e.$$
Therefore, 
$\set{\res(x)|x\in[\bl]\text{ for some }\bl\in\Lamb}\subseteq\Z/e\Z$. 

For $\bl\in\Lamb$ and $f\in\Z/e\Z$ put
$C_f(\bl)=\#\set{x\in[\bl]|\res(x)=f}$.  It is straightforward to
check that with these new conventions $\bl\sim_C\bmu$ if and only if
$C_f(\bl)=C_f(\bmu)$, for all $f\in\Z/e\Z$.

\subsection{Abacuses}
Abacuses first appeared in the work of Gordon
James~\cite{James:YoungD} and have since been used extensively in the
modular representation theory of the symmetric groups and related
algebras.  An \emph{$e$--abacus} is an abacus with $e$ vertical runners,
which are infinite in both directions. If $e$ is finite then we label
the runners $0,1,\dots,e-1$ from left to right and position $z\in\Z$
on the abacus is the bead position in row $x$ on runner~$y$, where $z=xe+y$ and $0 \leq y <e$. 
If $e=\infty$ then we label the runners $\ldots,-1,0,1,\ldots$ and
position $z$ on the abacus is the bead position in row $0$ on runner
$z$. 

Let $\bl\in\Lamb$ be a multipartition and recall that we have fixed a
sequence of integers $\c=(c_1,\dots,c_r)$. Fix $a$ with $1\le a\le r$
and, for $i\ge0$, define 
\[\beta^{a}_i = \lambda^{(a)}_i-i+c_a.\]
Then the \emph{$\beta$--numbers} (with charge $c_a$) for the partition
$\lambda^{(a)}$ are the integers $\beta_1^a,\beta_2^a,\ldots$ and we
define $B_a = \{\beta_1^a,\beta_2^a,\ldots \}$.  The \emph{$e$--abacus
display} of $\lambda^{(a)}$ is the $e$--abacus with a bead at
position~$\beta^{a}_i$, for $i\ge1$.  The $e$--abacus display of the
multipartition~$\bl$ is the ordered $r$--tuple of abacuses for the
partitions $\lambda^{(1)},\dots,\lambda^{(r)}$.

It is easy to check that a multipartition is uniquely determined by
its abacus display and that every abacus display corresponds to some
multipartition.

\begin{example}
Suppose that $e=3$, $r=3$ and $\c=(0,1,2)$. Let
$\bl=\((4,1,1),(2),(3,2,1)\)$. Then 
$$
B_1=\{3,-1,-2,-4,-5,\ldots\}, \quad
B_2=\{2,-1,-2,\ldots\},\quad
B_3=\{4,2,0,-2,-3,\ldots\}$$
and the abacus display for $\bl$ is given by
$$\begin{array}{c@{\hspace*{20mm}}c@{\hspace*{20mm}}c}
  \abacus{
    \vdots  & \vdots  & \vdots  & \\
      &   &   & \\
    - &   &   & \\
    - & - & - & \\
      & - & - & \\
    - & - & - & \\
     \, &\, &\, & \\
      \ncline{1,1}{7,1}
      \ncline{1,2}{7,2}
      \ncline{1,3}{7,3}
  } \qquad \qquad
   \abacus{
    \vdots  & \vdots  & \vdots  & \\
      &   &   & \\
      &   &   & \\
    -  & -  &  & \\
    - & - & - & \\
    - & - & - & \\
     \, &\, &\, & \\
      \ncline{1,1}{7,1}
      \ncline{1,2}{7,2}
      \ncline{1,3}{7,3}
  } \qquad \qquad
 \abacus{
    \vdots  & \vdots  & \vdots  & \\
      &   &   & \\
      &   & - & \\
      & - &   & \\
    -  &  & - & \\
    - & - & - & \\
     \, &\, &\, & \\
      \ncline{1,1}{7,1}
      \ncline{1,2}{7,2}
      \ncline{1,3}{7,3}
  }
\end{array}$$
\end{example}

Let $\lambda$ be a partition and suppose that
$B=\{\beta_1,\beta_2,\dots\}$ is the set of $\beta$--numbers for $\lambda$.
Then the $e$--abacus for $\lambda$ has beads at positions
$\beta_i$, for $i\ge0$. If $\beta_i+h\notin B$ then \emph{moving} the
bead at position $\beta_i$ to the \textit{right} $h$ positions gives a
new abacus display with beads at positions
$\{\beta_1,\beta_2,\dots,\beta_{i-1},\beta_i+h,\beta_{i+1},\dots\}$.
Similarly, if $\beta_i-h\notin B$ then \emph{moving} this bead $h$
positions to the \textit{left} creates a new abacus display with beads
at positions
$\{\beta_1,\beta_2,\dots,\beta_{i-1},\beta_i-h,\beta_{i+1},\dots\}$.
The conditions $\beta_i\pm h\notin B$ are needed to ensure that the abacus
display for $\lambda$ does not already have a bead at the new position.
Note that with these conventions moving a bead on runner~$0$ one position
to the left moves the bead to a position on runner $e-1$ in the preceding
row. Similarly, moving a bead on runner $e-1$ to the right moves a bead to
a position on runner~$0$ in the next row. We also talk of moving
beads in the abacus displays of multipartitions.

Recasting the above discussion in terms of the abacus we have the
following well--known result which goes back to Littlewood and James.
(Recall that we defined rim hooks in section~2.)

\begin{lemma}\label{wrapping}
Suppose that $\lambda$ is a partition. Then moving a bead to the
right $h$ positions from runner~$f$ to runner $f'$ corresponds to
wrapping an $h$--rim hook with foot residue $f$ onto $\lambda$.
Similarly, moving a bead~$h$ positions to the left, from runner $f$ to
runner $f'$ corresponds to unwrapping an $h$--rim hook from~$\lambda$
with foot residue $f$.
\end{lemma}

That increasing a $\beta$--number by $h$ corresponds to wrapping on an
$h$--rim hook is proved in~\cite[Lemma 5.26]{M:ULect}. The remaining
claim about residues follows easily from our definitions. As a
consequence we obtain the following.

\begin{corollary}\label{abacusprops}
Suppose that $\lambda$ is a partition and $f\in\Z/e\Z$, where $e<\infty$. 
Then
\begin{enumerate}
\item Moving a bead down one row on a runner
    corresponds to wrapping an $e$--rim hook onto $[\lambda]$. If this
    bead is on runner $f$ then the rim hook has foot residue~$f$.
\item Moving a bead up one row on a runner
    corresponds to unwrapping an $e$--rim hook from $[\lambda]$. If this
    bead is on runner $f$  the rim hook has foot residue~$f$.
\item Moving the lowest bead on runner $f$ down one row corresponds to
wrapping on an $e$--hook with foot residue $f$. Consequently, we can
add an $e$--hook with foot residue $f$ to any partition.  
\end{enumerate}
\end{corollary}

Suppose that $\lambda$ is a partition. The \emph{$e$--core} of
$\lambda$ is the partition $\bar\lambda$ whose $e$--abacus display is
obtained from the $e$--abacus display for $\lambda$ by moving all
beads as high as possible on their runners, that is, successively removing all $e$--hooks from the diagram of $\lambda$.  If $e=\infty$ then the
$e$--core of $\lambda$ is $\lambda$ itself. Define the \emph{$e$--weight} of the partition, $\ww(\lambda)$, to be the number of $e$--hooks that we remove in order to construct $\bar\lambda$.

\subsection{Jantzen equivalence}
In order to prove Theorem~\ref{blocks} we first simplify the formula
for $J_{\bl\bmu}$. Let $\bl$ be a multipartition and recall that if
$x\in[\bl]$ then $r^\bl_x\subseteq[\bl]$ is the associated rim hook.
To ease notation we let $h^\bl_x=|r^\bl_x|$ be the \emph{hook length} of ~$r^\bl_x$.

Recall that $\F$ is a field of characteristic~$p$. 
Define $\nup\map{\Z^\times}\N$ to be the map
$$\nup(h)=\begin{cases}
    p^k,&\text{if $p$ is finite},\\
    1,&\text{if $p=\infty$},
\end{cases}$$
where $k\ge0$ is maximal such that $p^k$ divides $h$. We caution the
reader that $\nup$ is not the standard $p$--adic valuation map.

If $\sigma=(\sigma_1,\sigma_2,\dots)$ is a partition let
$\sigma'=(\sigma_1',\sigma_2',\dots)$ be  its conjugate. Then
$\sigma_i'=c$ if $c$ is maximal such that $(c,i)\in[\sigma]$.  (So
$\sigma_i'$ is the length of column~$i$ of $[\sigma]$.) 
For any integer $h\in\Z$ let $[h]_t=\frac{t^h-1}{t-1}\in\F[t,t^{-1}]$.  

\begin{lemma}\label{valuation}
Suppose that $\bl$ and $\bmu$ are multipartitions of $n$ and that
$[\bl]{\setminus}r^\bl_x=[\bmu]{\setminus}r^\bmu_y$, for some
nodes $x=(i,j,a)\in[\bl]$ and $y=(k,l,b)\in[\bmu]$. 
Then $\nu_\pi\(\res_\O(f^\bl_x)-\res_\O(f^\bmu_y)\)\ne0$ if and only 
if~$\res(f^\bl_x)=\res(f^\bmu_y)$, in which case
$$\nu_\pi\(\res_\O(f^\bl_x)-\res_\O(f^\bmu_y)\)
   =\nup\(n(a-b)+j-\lambda^{(a)'}_i-l+\mu^{(b)'}_k\).$$
\end{lemma}

\begin{proof}
Let $i'=\lambda^{(a)'}_i$ and $k'=\mu^{(b)'}_k$ so that $f^\bl_x=(i',j,a)$
and $f^\bmu_y=(k',l,b)$. Then
\begin{align*}
\res_\O(f^\bl_x)-\res_\O(f^\bmu_y)
    &=q^{j-i'+c_a}t^{na+j-i'}-q^{l-k'+c_b}t^{nb+l-k'}\\
    &=q^{l-k'+c_b}t^{nb+l-k'}\(q^{j-i'-l+k'+c_a-c_b}
              t^{n(a-b)+j-i'-l+k'}-1\).
\end{align*}
Therefore, $\nu_\pi(\res_\O(x)-\res_\O(y))\ne0$ if and only if
$q^{j-i'-l+k'+c_a-c_b}=1$, which is if and only if
$\res(f^\bl_x)=q^{j-i'+c_a}=q^{l-k'+c_b}=\res(f^\bmu_y)$. 

Now suppose that
$\res(f^\bl_x)=\res(f^\bmu_y)$ and let $h=n(a-b)+j-i'-l+k'$. Note that
$h$ is non-zero because if $a=b$ then $h$ is the axial distance from $x$ to $y$. Then
 $$\nu_\pi\(\res_\O(f^\bl_x)-\res_\O(f^\bmu_y)\)
         =\nu_\pi(t^{n(a-b)+j-i'-l+k'}-1)
         =1+\nu_\pi([h]_t).$$
If $p=\infty$ then $(t-1)$ does not divide $[h]_t$, so that
$\nu_\pi(\res_\O(x)-\res_\O(y))=1=\nup(h)$. If $p$ is finite then
write $h=p^kh'$, where $p\nmid h'$. Then
$$[h]_t=[p^kh']_t=[p^k]_t[h']_{t^{p^k}}=(t-1)^{p^k-1}[h']_t^{p^k}.$$
Now, $t-1$ does not divide $[h']_t$ since $p\nmid h'$. Therefore, 
$\nu_\pi([h]_t)=\nup(h)-1$ and the result follows.
\end{proof}

We can now prove that (c)$\implies$(d) in Theorem~\ref{blocks}.

\begin{corollary}\label{c to d}
Suppose that $\bl\sim_J\bmu$, where $\bl,\bmu\in\Lamb$. Then $\bl\sim_C\bmu$.
\end{corollary}

\begin{proof}
Without loss of generality we may assume that $J_{\bl\bmu}\ne0$.
By Lemma~\ref{valuation} and Definition~\ref{Jlambdamu}, $J_{\bl\bmu}$ is
non--zero only if there exist nodes $x\in[\bl]$ and $y\in[\bmu]$ such
that $[\bl]{\setminus}r^\bl_x=[\bmu]{\setminus}r^\bmu_y$ and
$\res(f^\bl_x)=\res(f^\bmu_y)$. These two conditions imply that
$C_f(\bl)=C_f(\bmu)$, for all $f\in\Z/e\Z$, so that $\bl\sim_C\bmu$.
\end{proof}

Establishing the reverse implication in Theorem~\ref{blocks} takes
considerably more effort.  We start by explicitly describing the
Jantzen coefficients.

\begin{proposition}\label{J non-zero}
Let $\bl=(\lambda^{(1)},\dots,\lambda^{(r)})$ and
$\bmu=(\mu^{(1)},\dots,\mu^{(r)})$ be multipartitions in $\Lamb$.
\begin{enumerate}
\item Suppose that there exist integers $a< b$ such that
    $\lambda^{(c)}=\mu^{(c)}$, for $c\ne a,b$, and that 
    $\lambda^{(a)} \neq \mu^{(a)}$ and $\lambda^{(b)} \neq \mu^{(b)}$. Then $J_{\bl\bmu}\ne0$
    only if there exist nodes $x=(i,j,a)\in[\bl]$ and
    $y=(k,l,b)\in[\bmu]$ such that $\res(f^\bl_x)=\res(f^\bmu_y)$ and
    $[\bl]{\setminus}r^\bl_x=[\bmu]{\setminus}r^\bmu_y$. In this case
$$J_{\bl\bmu}=(-1)^{\leg(r^\bl_x)+\leg(r^\bmu_y)}
        \nup\(n(a-b)+j-\lambda^{(a)'}_i-l+\mu^{(b)'}_k\).$$
\item Suppose that $e$ is finite and for some integer $a$ we have
    $\lambda^{(c)}=\mu^{(c)}$, for $c\ne a$. Then $J_{\bl\bmu}\ne0$
    only if there exist nodes $x=(i,j,a),(i,m,a)\in[\bl]$ such that $m<j$, 
    $e\mid h^\bl_{(i,m,a)}$ and $\bmu$ is obtained by wrapping a rim hook 
    of length $h^\bl_x$ onto $\bl{\setminus}r^\bl_x$ with its hand node 
    in column~$m$. In this case
    $$J_{\bl\bmu}=\begin{cases}
	\space(-1)^{\leg(r^\bl_x)+\leg(r^\bmu_y)}
	\nup(h^\bl_{(i,m,a)}),&\text{if }e\nmid h^\bl_{(i,j,a)},\\[3mm]
	\space(-1)^{\leg(r^\bl_x)+\leg(r^\bmu_y)}
             \Big(\nup(h^\bl_{(i,m,a)})-\nup(h^\bl_{(i,j,a)})\Big),
	    &\text{if }e\mid h^\bl_{(i,j,a)},
	\end{cases}$$
    where the node $y\in[\bmu]$ is determined by
    $[\bmu]{\setminus}r^\bmu_y=[\bl]{\setminus}r^\bl_x$.
\item In all other cases, $J_{\bl\bmu}=0$.
\end{enumerate}
\end{proposition}

\begin{proof}
Suppose that $J_{\bl\bmu}\ne0$. Then $\bl\gdom\bmu$ by
Definition~\ref{Jlambdamu} and $\res(f^\bl_x)=\res(f^\bmu_y)$
by Lemma~\ref{valuation}. Furthermore, there exist nodes
$x=(i,j,a)\in[\bl]$ and $y=(k,l,b)\in[\bmu]$ such that
$[\bl]{\setminus}r^\bl_x=[\bmu]{\setminus}r^\bmu_y$. 
Consequently, $\lambda^{(c)}\ne\mu^{(c)}$ for at most two values of~$c$.
Therefore, since $\bl\gdom\bmu$, we may assume that we have integers $1\le a\le b\le r$ such
that $\lambda^{(c)}=\mu^{(c)}$, for $c\ne a,b$.

If $a\ne b$ then the nodes $x$ and $y$ are
uniquely determined because $r^\bl_x=[\lambda^{(a)}]{\setminus}[\mu^{(a)}]$ and
$r^\bmu_y=[\mu^{(b)}]{\setminus}[\lambda^{(b)}]$. Therefore, $\lambda^{(a)}\ne\mu^{(a)}$, 
$\lambda^{(b)}\ne\mu^{(b)}$ and we are in the situation considered in
part~(a). The formula for $J_{\bl\bmu}$ now follows directly from
Definition~\ref{Jlambdamu} and Lemma~\ref{valuation}.

Now assume that $a=b$. If $e=\infty$ then
$\res(f^\bl_x)=\res(f^\bmu_y)$ if and only if $x=y$ since
$h^\bl_x=h^\bmu_y$. This forces $\bl=\bmu$, which is not possible
since $\bl\gdom\bmu$. Hence, $e$ must be finite. 
By Lemma~\ref{wrapping}
the abacus display for $\mu^{(a)}$ is obtained from the abacus display
for $\lambda^{(a)}$ by moving one bead $h^\bl_x$ positions to the left
\textit{from} runner $\res(f^\bl_x)$, and other bead $h^\bl_x$
positions to the right \textit{to}  runner $\res(f^\bl_x)$.  

\Case{$e\nmid h^\bl_{(i,j,a)}$} By
Lemma~\ref{wrapping} and the remarks above, the beads on the abacus
displays of~$\lambda^{(a)}$ and $\mu^{(a)}$ are being moved between
different runners. Therefore, the nodes $x=(i,j,a)\in[\bl]$ and
$y=(k,l,a)\in[\bmu]$ are uniquely determined by the conditions
$\res(f^\bl_x)=\res(f^\bmu_y)$ and
$[\bl]{\setminus}r^\bl_x=[\bmu]{\setminus}r^\bmu_y$. Let
$m=\mu^{(a)}_k$. Then
$h^\bl_{(i,m,a)}=(j-\lambda^{(a)'}_i)-(l-\mu^{(a)'}_k)$ is the `axial
distance' from $f^\bl_x$ to $f^\bmu_y$, so that $e\mid h^\bl_{(i,m,a)}$. 
(In fact, $h^\bl_{(i,m,a)}$ is the axial distance between
the corresponding hand nodes, but this distance is, of course, the
same.  Note also that, since $\res(f^\bl_x)=\res(f^\bmu_y)$, we have that $e\mid h^\bl_{(i,m,a)}$.) Hence,
$J_{\bl\bmu}=(-1)^{\leg(r^\bl_x)+\leg(r^\bmu_y)}\nup\(h^\bl_{(i,m,a)}\)$
by Definition~\ref{Jlambdamu} and Lemma~\ref{valuation}. 

\Case{$e\mid h^\bl_{(i,j,a)}$} Since $h^\bl_x\equiv0\pmod
e$ unwrapping $r^\bl_x$ from $\bl$ and wrapping $r^\bmu_y$ back onto
$\bl{\setminus}r^\bl_x$ corresponds to moving one bead on
runner~$\res(f^\bl_x)$ up $\tfrac1eh^\bl_x$ rows and another bead on
runner~$\res(f^\bl_x)$ down $\frac1eh^\bl_x$ rows.  If in the abacus
display for $\bl$ these beads were moved from rows $r_1>r_2$ to rows~$r_1'$ 
and $r_2'$, respectively, then the abacus display for $\bmu$
can also be obtained from abacus display for~$\bl$ by moving the bead
in row $r_1$ to row $r_2'$ and moving the bead in row $r_2$ to row
$r_1'$. That is, there exist nodes $x'\neq x$ and $y'\neq y$ such that
we can obtain $\bmu$ by unwrapping $r^\bl_{x'}$ from $\bl$ and
wrapping $r^\bmu_{y'}$ back onto~$\bl{\setminus}r^\bl_{x'}$.
By Lemma~\ref{wrapping} there are no other ways of obtaining
$\bmu$ by unwrapping a rim hook from $\bl$ and wrapping it back on
again. Since $\bl\gdom\bmu$ we can choose the nodes $x=(i,j,a)$
and $y=(k,l,a)$ above so that $r_1>r_1'>r_2'>r_2$. Then
$x'=(i,m,a)$, where $m=\mu^{(a)}_k$, and $y'=(\lambda^{(a)'}_j,l,a)$.
Further, $\leg(r^\bl_x)+\leg(r^\bmu_y) = \lambda^{(a)'}_j-i + \mu^{(a)'}_l-k$ and $\leg(r^\bl_{x'})+\leg(r^\bmu_{y'}) = \lambda^{(a)'}_{m}-i + \mu^{(a)'}_l-\lambda^{(a)'}_j$.  By 
construction, $k=\lambda^{(a)'}_{m}+1$, so 
$\leg(r^\bl_{x})+\leg(r^\bmu_{y})$ and 
$\leg(r^\bl_{x'})+\leg(r^\bmu_{y'})$ have opposite parities. The axial
distance from $f^\bl_x$ to $f^\bmu_y$ is $h^\bl_{(i,m,a)}$ (where $e \mid h^\bl_{(i,m,a)}$ since $\res(f^\bl_x)=\res(f^\bmu_y)$) and the
axial distance from $f^\bl_{x'}$ to $f^\bmu_{y'}$ is
$h^\bl_{(i,j,a)}$. Therefore,
$$J_{\bl\bmu}=(-1)^{\leg(r^\bl_x)+\leg(r^\bmu_y)}
              \Big(\nup(h^\bl_{(i,m,a)})-\nup(h^\bl_{(i,j,a)})\Big)$$
as required.

We have now exhausted all of the cases where $J_{\bl\bmu}$ is
non--zero, so the Proposition is proved.
\end{proof}

\subsection{Residue equivalence}
We are now ready to start proving that $\bl\sim_J\bmu$ whenever
$\bl\sim_C\bmu$.

A rim hook of $\bl$ is \emph{vertical} if it is contained within a
single column of $[\bl]$.  

A partition $\lambda$ is an \emph{$(e,p)$--Carter partition} if it has the property that 
$$\nup(h^\lambda_{(i,m,1)})=\nup(h^\lambda_{(i,j,1)}), 
                   \text{\qquad for all }(i,m,1),(i,j,1)\in[\lambda].$$ 
These partitions arise because $\Delta(\lambda)$ is irreducible if and
only if $\lambda$ is $(e,p)$--irreducible. The
$(e,p)$--Carter partitions are described explicitly in
\cite[Theorem~5.45]{M:ULect}. For us the most important properties
of these partitions are that if $\lambda$ is an
$(e,p)$--Carter partition then:
\begin{itemize}
\item all of the $e$--hooks which can be unwrapped from
$\lambda$ when constructing its $e$--core $\bar{\lambda}$ are vertical;
\item $\nup$ is constant on the rows of $[\lambda]$; and 
\item $\bar{\lambda}'_i\equiv\bar{\lambda}'_{i-1}-1\pmod e$ whenever
$\lambda'_i\ne\bar{\lambda}'_i$. 
\end{itemize}

\begin{proposition} \label{Amoves}
Suppose that 
$\bl \in \Lamb$ and $1\leq a \leq r$.  Define \[\Lambda_a(\bl)=
\set{\bmu\in\Lamb|\bar{\mu^{(a)}}=\bar{\lambda^{(a)}}
\text{ and }  
\mu^{(c)}=\lambda^{(c)}
\text{ when }
c\ne a}.\]  
Then $\bl\sim_J\bmu$ for all $\bmu \in \Lambda_a(\bl)$.  
\end{proposition}

\begin{proof} 
Suppose that $\bmu \in \Lambda_a(\bl)$.  
If $e=\infty$ then $\bar{\lambda^{(a)}} = \bar{\mu^{(a)}}$ if and only if $\lambda^{(a)} = \mu^{(a)}$ so there is nothing to prove.
Assume then that $e$ is finite and let $w_a=\ww(\lambda^{(a)})$. If
$w_a=0$ then $\bar{\lambda^{(a)}}=\lambda^{(a)}$ so that $\bl=\bmu$ and
there is nothing to prove. So we can assume that $w_a>0$.

Let $\brho$ be the multipartition in $\Lambda_a(\bl)$ where $\rho^{(a)}$
is the partition obtained by wrapping $w_a$ vertical $e$--hooks
onto the first column of the $e$--core of $\lambda^{(a)}$. Then
$\bmu\gedom\brho$ for all $\bmu\in\Lambda_a(\bl)$. To prove the Lemma
it is enough to show that $\bmu\sim_J\brho$, for all
$\bmu\in\Lambda_a(\bl)$. By induction on dominance
we may assume that $\bnu\sim_J\brho$ whenever $\bnu\in\Lambda_a(\bl)$
and $\bmu\gdom\bnu$. 
If $J_{\bmu\bnu}\ne0$ for some $\bnu\in\Lambda_a(\bl)$
then $\bmu\sim_J\bnu$.  As $\bmu\gdom\bnu$, we have that
$\bnu\sim_J\brho$ by induction, so that $\bmu\sim_J\bnu\sim_J\brho$.

It remains to consider the case when $\bmu\gdom\brho$ and $J_{\bmu\bnu}=0$ for all
$\bnu\in\Lambda_a(\bl)$. By Lemma~\ref{J non-zero} (b),
$$\nup(h^\bmu_{(i,m,a)})=\nup(h^\bmu_{(i,j,a)}), 
                   \text{\qquad for all }(i,m,a),(i,j,a)\in[\bmu],$$
so that $\mu^{(a)}$ is an $(e,p)$--Carter partition. 
Since $w_a>0$ we can find a (unique) node $(i,j,a)\in[\bmu]$ such that
$h^\mu_{(i,j,a)}\equiv0\pmod e$ and $h^\mu_{(i',j',a)}\not\equiv0\pmod e$, 
for all $(i',j',a)\in[\bmu]$ with $(i',j')\ne(i,j)$, $i'\le i$ and
$j'\ge j$.  Let $\bnu$ be the multipartition obtained by unwrapping
$r^\bmu_{(i,j,a)}$ from $[\bmu]$ and wrapping it back on to the end of
the first row of $[\bmu]{\setminus}r^\bmu_{(i,j,a)}$. Similarly, let
$\bet$ be the multipartition obtained by unwrapping this same hook
from $\bmu$ and wrapping it back on to the end of the first column of
$[\bmu]{\setminus}r^\bmu_{(i,j,a)}$.   Therefore,
$J_{\bnu\bmu}\ne0$ and $J_{\bnu\bet}\ne0$, by Lemma~\ref{J non-zero} (b), 
so that $\bmu\sim_J\bnu\sim_J\bet$. Note that $\bmu\gdom\brho$ implies that $j>1$, so that $\bmu \gdom \bet$. Consequently,
$\bmu\sim_J\brho$ by induction.
\end{proof}

Recall that the $e$--cores of the partitions of $n$ completely determine the
blocks when $r=1$. We have the following imperfect generalization when
$r>1$.

\begin{definition}Suppose that
    $\bl=\(\lambda^{(1)},\dots,\lambda^{(r)}\)$ is a multipartition. 
    Then the \emph{$e$--multicore} of $\bl$ is the multipartition
    $\bar\bl=\(\bar\lambda^{(1)},\dots,\bar\lambda^{(r)}\)$.
We abuse notation and say that $\bl$ is a multicore if $\bl = \bar{\bl}$.
\end{definition}

By Corollary~\ref{abacusprops} (a), the $e$--multicore $\bar\bl$ of $\bl$
is obtained from $\bl$ by sequentially unwrapping all $e$--rim hooks
from the diagram of $\bl$, in any order. Note that if $e=\infty$ then
every multipartition is an $e$--multicore.

Mimicking the representation theory of the symmetric groups, we extend the definition of $\we$ to multipartitions by defining
$\we(\bl)$ to be the number of $e$--hooks that have to be unwrapped
from $\bl$ to construct $\bar{\bl}$. If $e$ is finite then 
$\we(\bl)=\frac1e(|\bl|-|\bar\bl|)$, whereas $\winfinity(\bl)=0$. Now define 
\[\We(\bl) = \max \{\we(\bmu) \mid \bmu \sim_C \bl \}.\]
Note that while $\We(\bl)$ is well defined, it is not immediately
clear how to compute it. 

\begin{lemma} \label{samecore}
Suppose that $\bl, \bmu \in \Lamb$ and that $\bar\bl=\bar\bmu$.
Then $\bl \sim_J \bmu$.
\end{lemma}

\begin{proof}We argue by induction on $d(\bl,\bmu)$, where 
$d(\bl,\bmu)=\frac1{e^2}\sum_{a=1}^r\(|\lambda^{(a)}|-|\mu^{(a)}|\)^2$. 
Note that $d(\bl,\bmu)$ is a non--negative integer because our assumption 
$\bar\bl=\bar\bmu$ implies that $|\lambda^{(a)}|\equiv|\mu^{(a)}|\pmod e$,
for $1\le a\le r$.

Suppose first that $d(\bl,\bmu)=0$. Then $|\la^{(a)}|=|\mu^{(a)}|$, for $1\le a\le r$. Define 
a sequence of multipartitions $\bnu_0=\bl, \bnu_1,\ldots,\bnu_r=\bmu$ by setting
\[\bnu^{(j)}_i = \begin{cases} 
\lambda^{(j)}, & i<j, \\
\mu^{(j)}, & i \geq j. 
\end{cases}\]
Then $\bnu_i \sim_J \bnu_{i+1}$ for $0 \leq i < r$, by Proposition~\ref{Amoves}, so that
$\bl\sim_J\bmu$ by transitivity.  

Now suppose that $d(\bl,\bmu)>0$. Since $\bar\bl=\bar\bmu$ and
$|\bl|=|\bmu|$, there exist integers $b$ and $c$ such that
$|\la^{(b)}|<|\mu^{(b)}|$ and $|\la^{(c)}|>|\mu^{(c)}|$. By
Corollary~\ref{abacusprops} it is possible to construct a new
multipartition~$\bnu$ by unwrapping an $e$--hook from~$\la^{(c)}$ and
wrapping it back on to $\la^{(b)}$ without changing the residue of the
foot node. Then $\bl\sim_J\bnu$ by Proposition~\ref{J non-zero} (and
Lemma~\ref{wrapping}). Moreover, $\bar\bnu=\bar\bl=\bar\bmu$ and
$d(\bnu,\bmu)<d(\bl,\bmu)$. Therefore, $\bnu\sim_J\bmu$ by induction,
so that $\bl\sim_J\bmu$ as required.
\end{proof}

We now need several results and definitions of Fayers from the
papers~\cite{Fayers:AKweight,Fayers:multicore}.  It should be noted
that there is a certain symbiosis between these two papers and the
present paper because Fayers wrote his papers believing that the
classification of the blocks of the Ariki--Koike algebras had already
been established. Fortunately, Fayers' results do not depend on the
block classification so when he discovered that there was a gap in the
previous proof of the classification he changed his papers so that
they now refer to `combinatorial blocks', or residue classes of
multipartitions. Thanks to the main result of this paper, the
`combinatorial blocks' studied by Fayers are indeed blocks.

\begin{definition}
\begin{enumerate}
\item 
(Fayers~\cite{Fayers:multicore}) 
Suppose that $\bl$ is a multicore and, if $e= \infty$, suppose further that 
the abacus display for $\lambda^{(a)}$ contains a bead in position
    $i$ but not in position $j$, while the abacus display for
    $\lambda^{(b)}$ contains a bead in position $j$ but not in
    position $i$.  Define
    $s^{ab}_{ij}(\bl)$ to be the multicore whose abacus display is
    obtained by moving a bead from runner $i$ to runner $j$ on the
    abacus for~$\lambda^{(a)}$ and moving a bead from runner $j$ to
    runner $i$ on the abacus for $\lambda^{(b)}$.
\item Suppose that $e$ is finite and let $\bl$ be a multipartition.
Define $t^a_{iw}(\bl)$ to be the multipartition whose abacus display
is obtained by moving the lowest bead on runner~$i$ of the abacus 
for $\lambda^{(a)}$ down~$w$ rows.
\end{enumerate}
\end{definition}

\begin{lemma} \label{stepeq}
Suppose that $\bl \sim_C \bmu$ and that $\bar{\bmu} = s^{ab}_{ij}(\bar{\bl})$.  Then $\bl \sim_J \bmu$.
\end{lemma}

\begin{proof}
Let $\bnu = t^a_{i \we(\bl)}(\bar{\bl})$ and $\brho=t^a_{j
\we(\bmu)}(\bar{\bmu})$. Then $\bl \sim_J \bnu$ and $\brho\sim_J \bmu$
by Lemma~\ref{samecore}. Furthermore, the multipartitions $\bnu$ and
$\brho$ satisfy the conditions of Proposition~\ref{J non-zero} (a), so
$\bl\sim_J\bnu\sim_J\brho\sim_J\bmu$ as required.
\end{proof}

\begin{definition}[\protect{Fayers~\cite[\S2.1]{Fayers:AKweight}}]
Suppose that $\bl$ is a multipartition. Then the \emph{$e$--weight} of
$\bl$ is the integer 
\[\Wt(\bl) = \sum_{j=1}^r C_{c_j}(\bl)-
          \frac{1}{2}\sum_{f\in\Z/e\Z}(C_f(\bl) - C_{f+1}(\bl))^2.\]
\end{definition}

Fayers~\cite{Fayers:AKweight} shows that $\Wt(\bl)\ge0$ for all
multipartitions $\bl$, and that $\Wt(\lambda)=\ww(\lambda)$ when
$r=1$;  that is, Fayers' definition of weight coincides with usual
definition of weight on the set of partitions. Further, if
$\bl\sim_C\bmu$ then $\Wt(\bl)=\Wt(\bmu)$, so the function
$\Wt(\cdot)$ is constant on the residue classes of~$\Lamb$. The
results of~\cite[Prop.~3.8]{Fayers:AKweight} show how to use the
abacus display of~$\bl$ to calculate $\Wt(\bl)$.  Combining this
method with Lemma~\ref{weightlink} below gives a way of computing
$\We(\bl)$ using the abacus display of~$\bl$. We leave the details to
the reader. 

Recall that a node $(i,j,a) \in [\bl]$ is \emph{removable} if
$[\bl]{\setminus}\{(i,j,a)\}$ is the diagram of some multipartition
in~$\Lambda_{r,n-1}^+$. Similarly, a node $(i,j,a) \notin [\bl]$
is \emph{addable} if $[\bl] \cup \{(i,j,a)\}$ is the diagram
of some multipartition in~$\Lambda_{r,n+1}^+$. The node
$x=(i,j,a)$ is an $f$--node if $\res(x)=f$.

Let $\bl$ be a multipartition. For $f\in\Z/e\Z$ and $a\in\{1,\ldots,r\}$, define
\[ \delta^a_f(\bl) = 
    \#\{\text{ removable $f$--nodes of } [\lambda^{(a)}]\,\}
      -\#\{\text{ addable $f$--nodes of } [\lambda^{(a)}]\,\}
\]
and set 
\[\delta_f(\bl)=\sum_{j=1}^r \delta^j_f(\bl).\]
The sequence $(\delta_f(\bl) \mid f\in\Z/e\Z)$ is the
\emph{hub} of $\bl$. The hub of $\bl$ can be read off the abacus
display of~$\bl$ using Lemma~\ref{wrapping}.

Observe that Corollary~\ref{abacusprops} implies that if $e$ is finite
then the hub is unchanged by wrapping $he$--hooks onto~$[\bl]$, for
$h\ge1$.  Furthermore, $\bl$ and $\bmu$ have the same hub if $\bmu =
s^{ab}_{ij}(\bl)$, for some $a,b,i,j$.

\begin{proposition}%
[\protect{Fayers~\cite[Proposition 3.2 and Lemma 3.3]{Fayers:AKweight}}]
\label{hub}
Suppose that $\bl$ is a multipartition of~$n$ and $\bmu$ is a
multipartition of~$m$. Then 
\begin{enumerate}
\item If $e<\infty$ and $\bl$ and $\bmu$ have the same hub then $m \equiv n \mod e$ and
\[\Wt(\bl) - \Wt(\bmu) = \frac{r(n-m)}{e} ;\]
\item If $n=m$ then $\bl \sim_C \bmu$ if and only if they have the same hub.   
\end{enumerate}
\noindent Consequently, if $\bmu$ is obtained from $\bl$ by wrapping
on an $e$--hook, then $\Wt(\bmu) = \Wt(\bl)+r$.
\end{proposition}

The next result will let us determine when $\We(\bl)=\we(\bl)$.

\begin{proposition}[\protect{Fayers~\cite[Theorem 3.1]{Fayers:multicore}}]\label{3.1}
Suppose that $\bl\in\Lamb$ is a multipartition. Then the following are equivalent.
\begin{enumerate}
\item $\bmu$ is a multicore whenever $\bmu\sim_C\bl$.
\item $\Wt(\bmu)\ge\Wt(\bl)$ whenever $\bmu$ and $\bl$ have the same hub.
\end{enumerate}
\end{proposition}

\begin{definition}
A multipartition $\bl$ is a \emph{reduced multicore} if it satisfies the
conditions of Proposition~\ref{3.1}.
\end{definition}

Not every multicore is reduced. If $\bl$ is a reduced multicore
then the block which contains $\Delta(\bl)$ is, in general, not
simple.  In contrast, when $r=1$ every core is a reduced multicore and
the block containing a core is always simple. If $\bl$ is an reduced
multicore then Fayers~\cite{Fayers:multicore} calls the set of
multipartitions $\set{\bmu|\bmu\sim_C\bl}$ a `core block'.

\begin{lemma} \label{weightlink}
Suppose that $\bl \in \Lambda_{n,r}^+$.  Then $\bar{\bl}$ is a
reduced multicore if and only if $\we(\bl) = \We(\bl)$.
\end{lemma}

\begin{proof}
Suppose $\we(\bl) \neq \We(\bl)$.  By definition, there exists a
multipartition $\bmu$ such that $\bmu \sim_C \bl$ and
$\we(\bmu)>\we(\bl)$.  Now $\bar{\bmu}$ and $\bar{\bl}$ have the same
hub, and by Proposition~\ref{hub}, $\Wt(\bar{\bmu}) < \Wt(\bar{\bl})$,
contradicting Condition (b) of Proposition~\ref{3.1}. Therefore, $\bar{\bl}$
is not a reduced multicore.

Now suppose that $\bar{\bl}$ is not a reduced multicore.
Then there exists a multipartition $\bmu$, which is not
a multicore, such that $\bmu \sim_C \bar{\bl}$. Let 
$\bnu=t^1_{0\we(\bl)}(\bmu)$.  Then $\bnu \sim_C \bl$ and 
$\we(\bnu)>\we(\bl)$. Hence, $\We(\bl)>\we(\bl)$.
\end{proof}

\begin{lemma}
[\protect{Fayers~\cite[Proof of Proposition 3.7 (1)]{Fayers:multicore}}]
\label{noncoremove}
Suppose that $\bl$ is a multicore which is not reduced.
Then there exists a sequence of multicores
$\bl_0=\bl,\bl_1,\dots,\bl_k=\bmu$ such that $\Wt(\bmu)< \Wt(\bl)$, and
$\bl_{m+1}=s^{a_m b_m}_{i_m j_m}(\bl_{m})$ and $\Wt(\bl_m) \leq
\Wt(\bl)$, for $0\le m<k$.
\end{lemma}

\begin{lemma}
[\protect{Fayers~\cite[Proof of Proposition 3.7 (2)]{Fayers:multicore}}]
\label{coremove}
Suppose that $\bl$ and $\bmu$ are reduced multicores and that
$\bl\sim_C\bmu$.  Then there exists a sequence of multicores
$\bl_0=\bl,\bl_1,\dots,\bl_k=\bmu$ such that 
$\bl_{m+1}=s^{a_mb_m}_{i_m j_m}(\bl_{m})$ and $\bl_{m+1} \sim_C \bl_m$, 
for $0 \leq m <k$.
\end{lemma}

We can now complete the proof of Theorem~\ref{blocks} when $q\ne1$ and
the parameters $Q_1,\dots,Q_r$ are non--zero. Consequently, this completes
the proofs of Theorem~A from the introduction.

\begin{theorem}\label{case1}
Suppose that $q\ne1$ and that the parameters $Q_1,\dots,Q_r$ are
non--zero. Let $\bl$ and $\bmu$ be multipartitions in
$\Lambda_{n,r}^+$. Then $\bl \sim_C \bmu$ if and only if $\bl \sim_J \bmu$.
\end{theorem}

\begin{proof} By Corollary~\ref{c to d} if $\bl\sim_J\bmu$ then
    $\bl\sim_C\bmu$. Suppose then that $\bl\sim_C\bmu$. To show that
    $\bl\sim_J\bmu$ it is sufficient to prove the following two
    statements.  Let $\bnu \in \Lamb$.
\begin{enumerate}
\item Suppose that $\we(\bnu) < \We(\bnu)$.  Then there exists $\bet\in\Lamb$ such that $\bet \sim_J \bnu$ and 
$\we(\bet) > \we(\bnu)$.
\item Suppose that $\bnu \sim_C \bet$ and that $\we(\bnu) = \We(\bnu) = \we(\bet)$.   Then $\bet \sim_J \bnu$.
\end{enumerate}
Suppose, as in (a), that $\we(\bnu) < \We(\bnu)$. Then $e$ is finite and by
Lemma~\ref{weightlink}, $\bar{\bnu}$ is not a reduced multicore. By
Lemma~\ref{noncoremove}, there exists a sequence of multicores
$\bnu_0=\bar{\bnu},\bnu_1,\dots,\bnu_k$ such that $\Wt(\bnu_k)<
\Wt(\bar{\bnu})$ and for $0 \le m<k$ we have $\bnu_{m+1}=s^{a_m b_m}_{i_m j_m}(\bnu_{m})$ and
$\Wt(\bnu_m) \leq \Wt(\bar{\bnu})$. For all $m$ with $0 \leq m \leq k$, we have that $\bnu_m$ and $\bar{\bnu}$ have the same hub, so Proposition~\ref{hub}
says that $|\bnu_m| \leq |\bar{\bnu}|$, that
$|\bar{\bnu}|\equiv|\bnu_m|\pmod e$ and that $|\bnu_k| < |\bar{\bnu}|$.  
Define $w_m=\we(\bnu) +\tfrac1e(|\bar{\bnu}|-|\bnu_m|)$ and set $\bet_m = t^1_{0
w_m}(\bnu_m)$ and $\bet=\bet_k$.  Then $\bet_m \sim_J \bet_{m+1}$, for $0 \leq m < k$, by Lemma~\ref{stepeq},
so that by Lemma~\ref{Amoves} and transitivity, $\bnu\sim_J\bet_0 \sim_J \bet_m=\bet$.  Moreover, 
$\we(\bet)= \we(\bnu)+\tfrac1e(|\bar{\bnu}|-|\bnu_k|)>\we(\bnu)$ as
required.

Now consider (b), that is, suppose that $\bnu \sim_C \bet$ and $\we(\bnu) = \We(\bnu) = \we(\bet)$.  By Lemma \ref{weightlink}, $\bar\bnu$ and $\bar\bet$ are reduced
multicores. Then, by Lemma~\ref{coremove},
there exist multicores  $\bnu_0=\bar{\bnu},\bnu_1,\dots,\bnu_k=\bar{\bet}$
such that $\bnu_{m+1}=s^{a_m b_m}_{i_m j_m}(\bnu_{m})$ and $\bnu_{m+1}
\sim_C \bnu_m$ for $0 \leq m <k$.  For $0 \leq m \leq k$, define $\bxi_m = t^1_{0
\we(\bnu)}(\bnu_m)$.  Then by Lemma~\ref{stepeq}, $\bxi_m \sim_J
\bxi_{m+1}$ and by Lemma~\ref{Amoves} and transitivity, $\bnu\sim_J\bxi_0\sim_J\bxi_k\sim_J\bet$ 
as required.
\end{proof}

\section{The blocks for algebras with exceptional parameters}
In this section we classify the blocks of the Ariki--Koike algebras
for the remaining cases from (\ref{cases}). That is, we assume
that the parameters satisfy one of the following four cases:
$$\begin{tabular}[t]{llp{60mm}}
    Case 2.& $r=1$ and $q=1$.&\\
    Case 3.& $r>1$, $q=1$ and $Q_1=\dots=Q_r=1$.\\
    Case 4.& $r>1$, $q=1$ and $Q_1=\dots=Q_r=0$.\\
    Case 5.& $r>1$, $q\ne1$ and $Q_1=\dots=Q_r=0$.
\end{tabular}$$
As in the previous section the basic strategy is to use the Jantzen
sum formula to analyze the combinatorics of the Jantzen coefficients.

We distinguish between cases 2 and 3 because the blocks differ
dramatically in these two cases. In fact, the blocks in case~2 behave
like the blocks when $q\ne1$ and the parameters $Q_1,\dots,Q_r$ are
non-zero. Quite surprisingly, the algebras $\H_{r,n}$ and $\Sch_{r,n}$
have only one block in cases~3--5.

In all cases the blocks of the algebras $\H_{r,n}$ and $\Sch_{r,n}$
are determined by Jantzen equivalence by
Proposition~\ref{Jantzen=blocks}. This section gives an explicit
description of when two multipartitions are Jantzen equivalent in 
cases 2--5 above.

\subsection{The blocks when $r=1$ and $q=1$}
Assume that we are in case~2 above and let $\H_n=\H_{1,n}$ and
$\Sch_n=\Sch_{1,n}$. In this case the Specht modules and Weyl
modules are indexed by partitions, rather than multipartitions, so we
write $\lambda$ in place of $\bl$, and so on.  The nodes in the
diagrams of partitions are all of the form $(i,j,1)$, for $i,j\ge1$,
so we drop the trailing $1$ from this notation and consider a node
to be an ordered pair $(i,j)$, so that
$[\lambda]=\set{(i,j)|1\le j\le\lambda_i}$.

As $q=1$ we have that $e=p$. Following section~3 define the
\emph{residue} of a node $x=(i,j)$ to be
$$\res(x)=(j-i)\pmod p.$$
Once again,
$\set{\res(x)|x\in[\bl]\text{ for some }\bl\in\Lamb}\subseteq\Z/p\Z$. 
For a partition $\lambda$ and $f\in\Z/p\Z$ put
$C_f(\bl)=\#\set{x\in[\bl]|\res(x)=f}$ and define $\lambda\sim_C\mu$
if $C_f(\lambda)=C_f(\mu)$, for all $f\in\Z/p\Z$. Then it is
well--known (and easy to prove using Corollary~\ref{abacusprops} (a))
that $\bl\sim_C\bmu$ if and only if $\lambda$ and $\mu$ have the same
$p$--core.

We can now prove Theorem~\ref{blocks} when $q=1$ and $r=1$.  To prove
this result we need to show that the Jantzen and residue equivalence
relations on the set of partitions coincide. We follow the argument of
the previous section.

The analogue of Lemma~\ref{valuation} in case~2 is as follows.

\begin{lemma}\label{valuation1}
Suppose that $\lambda$ and $\mu$ are partitions of $n$ and that
$[\lambda]{\setminus}r^\lambda_x=[\mu]{\setminus}r^\mu_y$, for some
nodes $x=(i,j)\in[\lambda]$ and $y=(k,l)\in[\mu]$. 
Then 
$$\nu_\pi\(\res_\O(f^\lambda_x)-\res_\O(f^\mu_y)\)
   =\nup\(j-\lambda'_j-l+\mu'_l\).$$
\end{lemma}

\begin{proof}
Let $i'=\lambda'_i$ and $k'=\mu'_k$ so that $f^\lambda_x=(i',j)$
and $f^\mu_y=(k',l)$. Then
$$ \res_\O(f^\lambda_x)-\res_\O(f^\mu_y)
    =t^{na+j-i'}-t^{na+l-k'}
    =t^{na+l-k'}(t^{j-i'-l+k'}-1).
$$
Mimicking the proof of Lemma~\ref{valuation}, let $h=j-i'-l+k'$. Then
 $$\nu_\pi\(\res_\O(f^\lambda_x)-\res_\O(f^\mu_y)\)
         =\nu_\pi(t^{j-i'-l+k'}-1)
         =1+\nu_\pi([h]_t).$$
Repeating the second half of the proof of Lemma~\ref{valuation} completes
the proof.
\end{proof}

The only difference between Lemma~\ref{valuation} and
Lemma~\ref{valuation1} is that now
$\nu_\pi\(\res_\O(f^\lambda_x)-\res_\O(f^\mu_y)\)$ is non--zero
whenever $[\lambda]{\setminus}r^\lambda_x=[\mu]{\setminus}r^\mu_y$;
that is, we no longer require that $\res(f^\lambda_x)=\res(f^\mu_y)$. 

\begin{proposition}\label{J non-zero1}
Let $\lambda$ and $\mu$ are partitions of $n$. Then $J_{\lambda\mu}$
is non--zero only if $p$ is finite and there exist nodes
$x=(i,j),(i,m)\in[\lambda]$ such that $m<j$, $p\mid h^\lambda_{(i,m)}$ 
and $\mu$ is obtained by wrapping a rim hook of length $h^\lambda_x$ onto
$\lambda{\setminus}r^\lambda_x$ with its highest node in column~$m$.
In this case 
$$J_{\lambda\mu}=\begin{cases}
    (-1)^{\leg(r^\lambda_x)+\leg(r^\mu_y)}\nup(h^\lambda_{(i,m)}),
        &\text{if }p\nmid h^\lambda_{(i,j)},\\[3mm]
    (-1)^{\leg(r^\lambda_x)+\leg(r^\mu_y)}
        \Big(\nup(h^\lambda_{(i,m)})-\nup(h^\lambda_{(i,j)})\Big),
	&\text{if }p\mid h^\lambda_{(i,j)},
\end{cases}$$
where the node $y\in[\mu]$ is determined by
$[\mu]{\setminus}r^\mu_y=[\lambda]{\setminus}r^\lambda_x$.
\end{proposition}

\begin{proof}
Suppose that $J_{\lambda\mu}\ne0$. Then $\lambda\gdom\mu$ by
Definition~\ref{Jlambdamu} and  there exist nodes
$x=(i,j)\in[\lambda]$ and $y=(k,l,b)\in[\mu]$ such that
$[\lambda]{\setminus}r^\lambda_x=[\mu]{\setminus}r^\mu_y$. 

\Case{$\res(f^\lambda_x)\ne\res(f^\mu_y)$} Unwrapping the rim hook
$r^\lambda_x$ from $\lambda$ moves a bead on the abacus for $\lambda$
from runner $\res(f^\lambda_x)$ to runner $r_1$, say, and wrapping the
rim hook $r^\mu_y$ back onto $\lambda{\setminus}r^\lambda_x$ moves a
bead from runner $r_2$ to runner $\res(f^\mu_y)$. Since
$\res(f^\lambda_x)\ne\res(f^\mu_y)$ we can also construct the
partition~$\mu$ from~$\lambda$ by moving a bead from runner
$\res(f^\lambda_x)$ to runner $r_2$ and then moving a bead from runner
$r_1$ to runner $\res(f^\mu_y)$. Comparing the abacus displays of
$\lambda$ and $\mu$, there are no other ways of obtaining $\mu$ from
$\lambda$ by moving a single rim hook.  As in the proof of Proposition \ref{J non-zero}, the sums of the leg lengths
for the two different ways of changing $\lambda$ into $\mu$ by moving
a rim hook have different parities, so their contributions to
$J_{\lambda\mu}$ cancel out. Hence, $J_{\lambda\mu}=0$ 
when~$\res(f^\lambda_x)\ne\res(f^\mu_y)$.

\Case{$\res(f^\lambda_x)=\res(f^\mu_y)$}The proof of
Proposition~\ref{J non-zero} in the case when $a=b$ can now be repeated
without change to complete the proof of the Proposition.
\end{proof}

\begin{corollary}\label{r=1 fini}
    Suppose that $\lambda$ and $\mu$ are partitions of
    $n$. Then $\lambda\sim_J\mu$ if and only if $\lambda\sim_C\mu$.
\end{corollary}

\begin{proof}By Proposition~\ref{J non-zero1}, $\lambda\sim_C\mu$ whenever 
    $\lambda\sim_J\mu$. The reverse implication follows by the argument of 
    Proposition~\ref{Amoves} since this proof only uses
    part~(b) of Proposition~\ref{J non-zero}, which is the same as 
    the statement of Proposition~\ref{J non-zero1}.
\end{proof}

\begin{remark}
Corollary~\ref{r=1 fini} completes the classification of the blocks of
the $q$--Schur algebras and the Hecke algebras of type~$A$; that is
when $r=1$.  Unfortunately, the classification of the blocks of the
$q$--Schur algebras given in \cite[Theorem~4.24]{JM:schaper} (and
reproduced in \cite[Theorem~5.47]{M:ULect}), contains a small error.
Fortunately, the classification of the blocks of the Hecke algebras of
type~$A$ given in \cite[Theorem~4.29]{JM:schaper} is correct --
indeed, when $r=1$ our proof is a streamlined version of this argument.
\end{remark}

\subsection{The blocks when $r>1$ and $q=1$ or $Q_1=\dots=Q_r=0$} We
now consider the blocks in the remaining cases, that is, when $r>1$
and either $q=1$ or $Q_1=\dots=Q_r=0$. In this case all simple modules
belong to the same block. We use the same strategy to prove
Theorem~\ref{blocks} in these cases as in the previous sections.  

Note that, in cases~3--5, $\res(x)=Q_a=Q_1$ for any node $x=(i,j,a)$.
Therefore, in these cases, $\Lamb$ forms a single residue class.
Hence, in order to prove Theorem~\ref{blocks}, we need to show that
any two multipartitions in $\Lamb$ are Jantzen equivalent.
Consequently, in cases~3--5 Theorem~\ref{blocks} asserts that the
algebras $\H_{r,n}$ and $\Sch_{r,n}$ have only one block. That is,
in cases~3--5, $\H_{r,n}$ and $\Sch_{r,n}$ are indecomposable algebras.

We adopt the same strategy that we used to prove Theorem~\ref{case1}.
To state the analogue of Lemma~\ref{valuation} set
$$\epsilon=\begin{cases}
       1,&\text{if $Q_1=\dots=Q_r=0$ (cases 4 and 5)},\\
       0,&\text{otherwise.}\\
\end{cases}$$

\begin{lemma}\label{valuation0}
Suppose that $\bl$ and $\bmu$ are multipartitions of $n$ and that
$[\bl]{\setminus}r^\bl_x=[\bmu]{\setminus}r^\bmu_y$, for some
nodes $x=(i,j,a)\in[\bl]$ and $y=(k,l,b)\in[\bmu]$. 
Then 
$$\nu_\pi\(\res_\O(f^\bl_x)-\res_\O(f^\bmu_y)\)=
    \nup\(n(a-b)+j-\bl^{(a)'}_i-l+\bmu^{(b)'}_k\)+\epsilon
$$
\end{lemma}

The proof of Lemma~\ref{valuation0} is similar to proofs of
Lemma~\ref{valuation} and Lemma~\ref{valuation1}, so we leave the
details to the reader. Note, in particular, that if $a\ne b$ then
$\nu_\pi\(\res_\O(f^\bl_x)-\res_\O(f^\bmu_y)\)$ is always non--zero
since $\nup(h)\ge0$, for all $h\in\Z\setminus\{0\}$. This crucial
difference leads to $J_{\bl\bmu}$ being non--zero whenever there exist
nodes $x=(i,j,a)\in[\bl]$ and $y=(k,l,b)\in[\bmu]$ with $a<b$ and
$[\bl]{\setminus}r^\bl_x=[\bmu]{\setminus}r^\bmu_y$. More explicitly,
we have the following analogue of Propositions~\ref{J non-zero}
and~\ref{J non-zero1}. Again, we leave details to the reader.

\begin{proposition}\label{J non-zero0}
Let $\bl=(\bl^{(1)},\dots,\bl^{(r)})$ and
$\bmu=(\bmu^{(1)},\dots,\bmu^{(r)})$ be multipartitions in $\Lamb$.
\begin{enumerate}
\item Suppose that there exist integers $a \neq b$ such that
    $\bl^{(c)}=\bmu^{(c)}$, for $c\ne a,b$. Then $J_{\bl\bmu}\ne0$
    only if $a<b$ and there exist nodes $x=(i,j,a)\in[\bl]$ and
    $y=(k,l,b)\in[\bmu]$ such that 
    $[\bl]{\setminus}r^\bl_x=[\bmu]{\setminus}r^\bmu_y$. In this case
$$J_{\bl\bmu}=(-1)^{\leg(r^\bl_x)+\leg(r^\bmu_y)}\Big(
        \nup\(n(a-b)+j-\bl^{(a)'}_i-l+\bmu^{(b)'}_k\)+\epsilon\Big).$$
\item Suppose that $e$ is finite and for some integer $a$ we have
    $\bl^{(c)}=\bmu^{(c)}$, for $c\ne a$. Then $J_{\bl\bmu}\ne0$
    only if there exist nodes $x=(i,j,a),(i,m,a)\in[\bl]$ such that
    $m<j$, $e\mid h^\lambda_{i,m,a)}$ and $\bmu$ is obtained by wrapping
    a rim hook of length $h^\bl_x$ onto $\bl{\setminus}r^\bl_x$ with its
    highest node in column~$m$. In this case
    $$J_{\bl\bmu}=\begin{cases}
	\space(-1)^{\leg(r^\bl_x)+\leg(r^\bmu_y)}
	\Big(\nup(h^\bl_{(i,m,a)})+\epsilon\Big),
	    &\text{if }e\nmid h^\bl_{(i,j,a)},\\[3mm]
	\space(-1)^{\leg(r^\bl_x)+\leg(r^\bmu_y)}
             \Big(\nup(h^\bl_{(i,m,a)})-\nup(h^\bl_{(i,j,a)})\Big),
	    &\text{if }e\mid h^\bl_{(i,j,a)},
	\end{cases}$$
    where $y\in[\bmu]$ is determined by
    $[\bmu]{\setminus}r^\bmu_y=[\bl]{\setminus}r^\bl_x$.
\item In all other cases, $J_{\bl\bmu}=0$.
\end{enumerate}
\end{proposition}

We can now complete the proof of Theorem~\ref{blocks}.

\begin{proof}[Proof of Theorem~\ref{blocks} for cases 3--5]
Let $\bl=(\lambda^{(1)},\dots,\lambda^{(r)})$ be a multipartition
of~$n$ and fix  integers $a\ne b$ with $\lambda^{(a)}\ne(0)$ and
$1\le a,b\le r$. Let $\bmu$ be any multipartition that can be obtained
by unwrapping a rim hook from $[\lambda^{(a)}]$ and wrapping it back
on to component~$b$ of $\lambda$. Then $\bl\sim_J\bmu$ by
Proposition~\ref{J non-zero0}(a). In particular, note that
$\bl\sim_J\bmu$ if $\bmu$ is obtained from $\bl$ by moving a removable
node from $\lambda^{(a)}$ to~$\lambda^{(b)}$. Consequently, by moving
the nodes in $[\bl]$ to the right, one by one, we see that $\bl$ is
Jantzen equivalent to a multipartition $\bmu$, where
$\bmu=((0),\dots,(0),\mu^{(r)})$ for some partition $\mu^{(r)}$.
Similarly, moving nodes in~$\bmu$ to the left, one by one, now shows
that $\bl\sim_J\bmu\sim_J((n),(0),\dots,(0))$.  Hence, every
multipartition in~$\Lamb$ is Jantzen equivalent
to~$((n),(0),\dots,(0))$. This shows that there is only one block in
cases~3, 4 and~5, so the Theorem follows.
\end{proof}

\section*{Acknowledgment}
We thank the referee for their detailed report and for
suggesting a number of improvements.


\end{document}